\documentclass[ejs]{imsart}

\usepackage{amsthm,amsmath,amssymb,latexsym,amsfonts,imsart}
\usepackage{comment}
\RequirePackage[colorlinks,citecolor=blue,urlcolor=blue]{hyperref}

\RequirePackage{hypernat}

\doi{10.1214/12-EJS693}
\pubyear{2012}
\volume{6}
\firstpage{861}
\lastpage{907}
\arxiv{0912.4489}

 \startlocaldefs
\renewcommand{\(}{$\,}
\renewcommand{\)}{\,$}

\newcommand{\cc}[1]{\mathcal{#1}}
\newcommand{\bb}[1]{\boldsymbol{#1}}

\renewcommand{\bar}[1]{\overline{#1}}
\renewcommand{\hat}[1]{\widehat{#1}}
\renewcommand{\tilde}[1]{\widetilde{#1}}
\newcommand{\nn}{\nonumber \\}
\numberwithin{equation}{section}
\numberwithin{figure}{section}
\newcounter{example}[section]
\numberwithin{example}{section}
\newcounter{remark}[section]
\numberwithin{remark}{section}
\newtheorem{theorem}{Theorem}[section]
\newtheorem{proposition}[theorem]{Proposition}
\newtheorem{definition}[theorem]{Definition}
\newtheorem{lemma}[theorem]{Lemma}
\newtheorem{corollary}[theorem]{Corollary}
\newtheorem{exmp}[example]{Example}
\newtheorem{rmrk}[remark]{Remark}
\newenvironment{example}{\begin{exmp}\rm}{\end{exmp}}
\newenvironment{remark}{\begin{rmrk}\rm}{\end{rmrk}}
\newcommand{\argmax}{\operatornamewithlimits{argmax}}

\def\argmin{\operatornamewithlimits{argmin}}

\def\Var{\operatorname{Var}}
\def\eqdef{\stackrel{\operatorname{def}}{=}}
\def\eqdistr{\stackrel{\operatorname{d}}{=}}

\def\E{I\!\!E}
\def\P{I\!\!P}
\def\R{I\!\!R}
\def\T{\top}
\def\kappa{\varkappa}
\newcommand{\diag}{\operatorname{diag}}
\newcommand{\I}{\mathbb{I}}
\def\ind{\mathbb{I}}

\newcommand{\s}{\sigma}
\newcommand{\bbpf}[1]{\tta^{*}_{#1}}
\newcommand{\fle}[2]{\tilde{f}_{#1}{(#2)}}
\newcommand{\flej}[3]{\tilde{f}^{(#1)}_{#2}{(#3)}}
\newcommand{\mle}[2]{\tilde{\theta}_{#1}^{(#2)}}
\newcommand{\mmle}[1]{\tilde{\tta}_{#1}}

\def\aadapest{\hat{\tta}}
\def\adapest{\hat{\theta}}
\def\adaplpest{\hat{f}}
\def\adapind{\hat{k}}
\def\dd{\mathrm{d}}

\def\LL{\operatorname{L}}

\def\EE{\mathbb{E}}
\def\PPi{\bb{\Pi}}
\def\RR{\mathbb{R}}
\def\RRn{\mathbb{R}^{n}}
\def\RRp{\mathbb{R}^{p}}
\def\RRd{\mathbb{R}^{d}}
\def\YY{\bb{Y}}
\def\Yi{Y_i}
\def\Xi{X_i}
\def\ff{\bb{f}}
\def\ffi{f(\Xi)}
\def\fta{f_{\bb\theta}}

\newcommand{\ta}[1]{\theta^{(#1)}}
\newcommand{\W}[1]{\mathbf{W}_{#1}}
\newcommand{\w}[2]{w_{{#1},{#2}}}
\newcommand{\B}[1]{\mathbf{B}_{#1}}
\def\KL{\mathbb{K}\mathbb{L}}
\def\M{\mathbf{M}}
\def\S{\mathbf{S}}
\def\Lam{\mathbf{\Lambda}}
\def\tta{\bb\theta}
\def\DD{\mathbf{D}}
\def\TTa{\mathbf{\Theta}}
\def\SSigma{\mathbf{\Sigma}}
\def\VV{\mathbf{V}}
\def\eps{\varepsilon}
\def\eeps{\bb\varepsilon}

\def\epsi{\varepsilon_{i}}

\newcommand{\norm}[2]{\cc{N} \left( {#1},{#2} \right)}
\def\p{p}

\def\cov{\operatorname{cov}}

\def\supp{\operatorname{supp}}
\def\rank{\operatorname{rank}}
\def\dim{\operatorname{dim}}
\def\tr{\operatorname{tr}}
\def\vec{\operatorname{vec}}
\def\PPsi{\bb{\Psi}}

\def\Psii{\Psi_{i}}
\def\z{\mathfrak{z}}
\endlocaldefs

\begin{document}

\begin{frontmatter}

\title{Spatial adaptation in \\heteroscedastic regression: \\Propagation approach\thanksref{t1}
}
 \thankstext{t1}{Funding of the DFG FOR 916 and ANR-07-BLAN-0234 is acknowledged }
\runtitle{Adaptation in heteroscedastic regression}

\author{\fnms{Nora} \snm{Serdyukova}\ead[label=e1]{nora.serdyukova@mathematik.uni-goettingen.de}
}
\thankstext{t2}{The author wishes to thank the Associate Editor and unknown referee for many fruitful questions and comments that greatly improved the paper, as well as her supervisor Professor Vladimir Spokoiny for introduction to the astonishing world of adaptive estimation.}
 \address{Institute for Mathematical
Stochastics, Georg August Universit\"at G\"ottingen, Goldschmidtstr. 7, 37077 G\"ottingen GERMANY\\ \printead{e1}
}


\runauthor{Nora Serdyukova}
\maketitle
\begin{abstract}
The paper concerns the problem of pointwise adaptive estimation in regression when the noise is heteroscedastic and incorrectly known.
The use of the local approximation method, which includes the local polynomial smoothing as a particular case, leads to a finite family of estimators corresponding to different degrees of smoothing. Data-driven choice of localization degree in this case can be understood as the problem of selection from this family. This task can be performed by a suggested in Katkovnik and Spokoiny (2008) FLL technique based on Lepski's method. An important issue with this type of procedures -- the choice of certain tuning parameters -- was addressed in Spokoiny and Vial (2009). The authors called their approach to the parameter calibration ``propagation''. In the present paper the propagation approach is developed and justified for the heteroscedastic case in presence of the noise misspecification. Our analysis shows that the adaptive procedure allows a misspecification of the covariance matrix with a relative error of order \(  (\log n)^{-1} \), where \( n \) is the sample size.
\end{abstract}




\tableofcontents
\end{frontmatter}

\section{Introduction}

Consider a regression model
\begin{equation}\label{true model}
    \YY = \ff + \Sigma_0^{1/2} \eeps,\;\;\;\;\eeps \sim \norm{0}{I_n}
\end{equation}
with response vector \( \YY \in \RRn \) and unknown diagonal covariance matrix \( \Sigma_0 = \)
\(  \diag(\sigma_{0,1}^2, \ldots ,\sigma_{0,n}^2 ) \).
Let \( \cc X \) be a Borel subset of \( \RRn \) and \( \Xi \) be fixed elements of \( \cc X \).
Denote by \( f:\cc X \to \RR \) the unknown regression function, then with \( \ff = (f(X_1), \ldots, f(X_n)  )^{\T} \) model \eqref{true model} can be written as
\begin{equation}\label{true model par coordonees}
     \Yi =\ffi + \sigma_{0,i} \,\epsi , \;\;\; i=1, \ldots , n.
\end{equation}
Given a point \( x \in \cc X\), the target of estimation is the value of \( f(x) \). The idea is to replace model \eqref{true model par coordonees} by a local parametric model
\begin{equation}\label{locparam model general}
    y_i = f_{\tta}(\Xi) + \sigma_{i} \,\epsilon_i , \;\;\; i : \Xi \in U_h(x),
\end{equation}
where \( \s_i > 0 \) are known, \( U_h(x) \eqdef \{ t : \| t-x \| \le h/2 \} \) and \( \tta \in \Theta \subset \RRp \) is an unknown parameter to be estimated. Denote by \( \PPsi = (\Psi_1 , \ldots , \Psi_n) \) a \( p \times n \) design matrix. In the considered set-up the covariance matrix \( \Sigma_0 \) is not known exactly and \( \Sigma= \diag(\sigma_1^2, \ldots , \sigma_n^2 ) \) stands for the available covariance matrix. Then the approximate model used instead of the true one reads as follows:
\begin{equation}\label{PA}
    \bb y = \PPsi^{\T}\tta + \Sigma^{1/2} \eeps.
\end{equation}
 Employing inside of \( U_h(x) \) one of the well-developed parametric methods we can estimate \( \tta \) by \( \tilde \tta (y_1, \ldots, y_d; x) \), and then use the estimator \(  f_{\tilde \tta (Y_1, \ldots, Y_d)}(x)  \) based on the observations from the ``true'' model \eqref{true model par coordonees} for estimation of \( f(x) \). Therefore we have to choose the local model
 (correspondingly, the collection of estimators of \( f_{\tta}(\cdot), \tta \in \Theta \)) and the appropriate degree of locality \( h \). This method of local approximation originated from \cite{Stone77}, \cite{Cleveland79}, \cite{Katk1979}, \cite{Stone80}, \cite{Katk1983}, \cite{Katk1985}, \cite{Tsybakov86}, \cite{KatkEA2006}.
In what follows we consider approximation by local linear models of the type:
\begin{equation}\label{locparam model linear}
    y_i = \Psii \tta + \sigma_{i} \,\epsilon_i , \;\;\; i : \Xi \in U_h(x),
\end{equation}
where \(\Psii= \Psi(\Xi) = (\psi_1(\Xi-x), \ldots, \psi_p(\Xi-x))^{\T}\) is a vector of basis functions
\( \{ \psi_j(\cdot) \} \) which already are fixed. Thus the model is misspecified in two places: in the form of the regression function and in the error distribution. The main issue then is to choose the appropriate bandwidth
\( h \) such that the estimator
 \begin{equation}\label{form of the estimator}
      f_{\tilde \tta_h} (x) \eqdef  \sum_{j=1}^p \tilde \theta^{(j)}_{h}(x) \psi_j(0)
 \end{equation}
built on the base of localized data would provide a relevant estimator for \( f(x) \). For this purposes the bandwidths selection should be done in a data-driven way, and this problem can be formulated as adaptive selection from the finite family \( \{ f_{\tilde \tta_h} (x) \}_{h>0} \). Notice also that the coefficients
\( \theta^{(j)}(x) \) as well as their estimators depend on \( x \) and should be calculated for every particular point of interest \( x \). On the other side the localization reduces influence of the choice of the functions
\( \{ \psi_j(\cdot) \} \) allowing to use simple collections.

 \par The proposed approach includes the important class of polynomial regressions, see
 \cite{Fan and Gijbels book}, \cite{KatkEA2006}, \cite{Loader}, \cite{Tsybakov}. For example in the univariate case \( x \in \R \), due to
the Taylor theorem, the approximation of the unknown function \( f(t) \) for \( t \) close to \( x \) can be written in the following form: \( \fta(t)= \ta{0} + \ta{1} (t-x) + \cdots + \ta{\p-1} (t-x)^{\p-1}/(p-1)! \) with the parameter \( \tta =(\ta{0},\,\ta{1}, \ldots , \ta{\p-1})^{\T} \) corresponding to the values of \( f(\cdot) \) and its derivatives at the point \( x \), if they exist. The design matrix~\( \PPsi \) then consists of the columns \begin{equation*}
      \Psii= \left(1,\, \Xi -x, \ldots, (\Xi -x)^{\p-1}/(p-1)!\right)^{\T} \;, \;\; i=1,\ldots ,n,
\end{equation*}
 and corresponds to the well
known polynomial smoothing. If the regression function is sufficiently smooth, then, up to a reminder term, for any \( t \) close to \( x \), \( f(t) \approx  \fta(t) \) and the estimator of \( f(x) \) at \( x \) is given by \( \fle{}{x} = f_{\tilde{\tta} (x)} (x) = \tilde{\theta}^{(0)} \). More details on the local polynomial estimation can be found, for instance, in \cite{Fan and Gijbels book}, \cite{Loader} or \cite{Tsybakov}. The local constant fit at a given point \( x \in \R \) is covered as well with \( p=1 \). In this case the ``design'' matrix is a row \(   \PPsi = (1, \ldots ,1) \) and \( \fta(\Xi) = \PPsi_i^{\T} \tta = \theta^{(0)} = \fta(x)  \), \( i=1, \ldots, n\). This type of approximation in our set-up with known constant noise is treated in \cite{KatkSpok} and \cite{SV}.

\par Nonparametric estimation in heteroscedastic regression under the \( L_2 \) losses was studied in \cite{EfroimovichPinsker96}, \cite{Efroimovich2007} and series of papers \cite{GaltchoukPergam2009}, \cite{GaltchoukPergam2010a}, \cite{GaltchoukPergam2010b}. One should mention very interesting paper \cite{DalalyanSalmon} on aggregation estimation under empirical losses in heteroscedastic Gaussian regression.
For estimation of the mean with \( L_2 \)-risk in Gaussian homoscedastic model with unknown variance the penalties allowing to deal with the complexity of such a collection of models were proposed in \cite{BaraudGiraudHuet}. However the problem of ``local model selection'' addressed in the present paper is quite different to the model selection in the sense of \cite{BirgeMassart} and \cite{Massart} related to estimation with global risk. In this set-up an amazing progress is achieved for the model selection in heteroscedastic not necessary Gaussian regression model in \cite{Arlot}, \cite{ArlotMassart}, \cite{Saumard}. The minimax pointwise estimation in heteroscedastic regression is in focus of \cite{Brua}.

\section{Estimation procedure}
\subsection{Local parametric estimation}
\par Using the conceptual framework given in the introduction we choose the maximum likelihood estimation as a parametric method used inside of a smoothing window. Let us briefly recall the idea of the local likelihood method dating back to \cite{Brillinger} and \cite{TibshiraniHastie}.

\par If the response variables \( \Yi \) are independent and have a density \( v(y,s(\Xi)) \), then the joint log-density of the sample is given by \( \LL(s) = \sum_{i=1}^n \log v(\Yi, s(\Xi)) \) leading to the ``global'' maximum likelihood estimation. Let as before \( \fta (\cdot) \) be a function entirely described by a vector \( \tta \in \TTa \subset \RRp \). The local likelihood model does not assume that \( s(\Xi) = \fta (\Xi)\), but one fits the ``parametric'' model locally within the smoothing window described by weights \( \cc{W} (x) = \{ w_{i}(x) \}_{i=1}^n \). The \emph{local log-likelihood} is defined as
\begin{equation}\label{def_local_log-likelihood}
    \LL(\cc W,\tta) = \sum_{i=1}^n \log v(\Yi, \fta (\Xi)) w_{i}(x).
\end{equation}
The local likelihood estimator \( \mmle{ }(x) \) is a maximizer of this weighted sum, \( \mmle{ }(x) = \argmax_{\tta}\LL(\cc W,\tta)  \). It is worth pointing out that in spite of the term ``local likelihood'' seems to be standard, see \cite{Loader} for example, if the weights \( w_{i} \) are allowed to take values different from zero and one, the quantity defined by \eqref{def_local_log-likelihood} is not a log-likelihood in the probabilistic sense even if the data indeed locally follows the parametric model with \( v(\Yi, \fta (\Xi)) \) for all \( i: \) \( w_{i}(x)>0 \). However, the local, or more correctly, weighted log-likelihood inherits most of useful properties from its ``global'' counterpart, c.f. Proposition \ref{Wilks theorem}. And -- what is of particular importance -- the true value of the parameter \( \tta \) maximizes the expectation of \eqref{def_local_log-likelihood}, see \cite{Loader} p.72. This property in more general set-up leads to the minimum contrast \( -\LL(\cc W,\tta) \) estimation.

\par Leaving the computational aspects aside, the key issue of this method is a proper choice of the largest smoothing window where the parametric fit \( \fta \) is still adequate. Putting differently, if we consider a finite collection of smoothing windows and corresponding (quasi) MLE's, the target is a data-driven selection from this family. In what follows we explore this approach.

\par Fix a point \( x \in \RR^d \) as a center of localization and basis \( \{ \psi_j \} \). Denote by
 \begin{equation*}
     \Psii = \Psi(\Xi -x) = (\psi_1(\Xi -x) , \ldots, \psi_p(\Xi -x))^{\T},\;\;\;i=1, \ldots, n,
\end{equation*}
the vectors of basis functions centered at \( x \).
For the next nonparametric ``selection'' step we need a sequence of nested windows. Let for every \( x \) a finite sequence of scales \( \cc{W}_{k}(x) \), \( k=1, \ldots, K \), be given by matrices
\begin{equation*}
     \cc{W}_{k}(x) = \diag(\w{k}{1}(x), \ldots , \w{k}{n}(x)),
\end{equation*}
where the weights
  \( \w{k}{i}(x)\in [0,1]\) can be understood, for instance, as smoothing kernels \( \w{k}{i}(x) = W ((\Xi -x) h_k^{-1}) \). A particular localizing function \( w_{(\cdot , \cdot )} (x) \) is assumed to be fixed; the aim is to choose on the base of available data an index \( k \) of an ``optimal'' scale. To simplify the notation we sometimes suppress the dependence on the reference point \( x \). Denote by
  \begin{equation}\label{def_Wk}
      \W{k} \eqdef  \Sigma^{-1/2} \cc{W}_{k}  \Sigma^{-1/2}  =
        \diag \left( \frac{\w{k}{1}}{\sigma_1^2} , \ldots ,
        \frac{\w{k}{n}}{\sigma_n^2}\right),\;\;\;k=1, \ldots, K.
  \end{equation}
Let \( \Theta \) be a compact subset of \( \RRp \). Inside of any window given by \( \cc{W}_{k} \),
\( k= 1, \ldots K \), according to~\eqref{def_local_log-likelihood} for each \( k \) we calculate the (quasi) MLE \( \mmle{k} =  \mmle{k}(x) =
(\mle{k}{0}(x), \ldots , \mle{k}{\,\p-1}(x))^{\T} \) of \( \tta \):
\begin{equation}\label{MLE}
    \mmle{k}  \eqdef \argmax_{\tta \in \Theta} \LL(\W{k},\tta)  ,
\end{equation}
where \( \LL(\W{k},\tta) \) is the weighted log-likelihood corresponding to the joint distribution of independent sample with \( \Yi \sim \norm{\Psii^{\T} \tta}{\s_i^{2}}\):
\begin{eqnarray}\label{log-likelihood}
      \LL(\W{k},\tta)
      &=& -\frac{1}{2}
                \sum_{i=1}^{n} |\Yi -\Psii^{\T} \tta|^2 \frac{\w{k}{i}}{ \sigma_i^2} +R \\ \nonumber
            &=&
            -\frac{1}{2 }
                ( \YY -\PPsi^{\T}\tta )^{\T} \W{k}  ( \YY -\PPsi^{\T}\tta ) + R.
\end{eqnarray}
Here \( R \) stands for the terms independent of \( \tta \).
If the \( p \times p \) matrix \( \B{k} = \B{k}(x) \) given by
\begin{equation}\label{B}
     \B{k} \eqdef   \PPsi  \W{k}\PPsi^{\T} =
        \sum_{i=1}^{n} \Psii \Psii^{\T} \frac{\w{k}{i}}{ \sigma_i^2}
\end{equation}
is positive definite at the point \( x \), \( \B{k}(x) \succ 0 \), then \( \mmle{k} = \mmle{k}(x) \) given by
\begin{eqnarray}\label{MLE_formula}
  \mmle{k}     &=& \B{k}^{-1} \PPsi \W{k} \YY
                = \B{k}^{-1} \sum_{i=1}^{n}  \Psii \Yi \frac{\w{k}{i}}{ \sigma_i^2}
\end{eqnarray}
is a linear estimator. Recall that in the case of polynomial basis \( \{ t^{q} \}_{q=0}^{p-1} \) for every fixed \( k \) the first coordinate of \( \mmle{k}(x) \) is the local polynomial estimator for the value of \( f(x) \).

In what follows we assume that \( n>p \) and \( \det \B{k}(x) > 0 \) for any \( k=1, \ldots , K \). One needs to keep in mind that for example, if \( \w{1}{\cdot} = W((\cdot - x)h_1^{-1})\) is a finitely supported kernel function, one can always find a bandwidth \( h_1 \) so small that the matrix
\( \B{1}(x) \) is degenerated. This implies that the smallest value of \( h_1 \) should be chosen in order to guarantee \( \B{1}(x) \succ 0 \). More precisely we assume the following:
\begin{description}
\item[\( \bb{(A1)} \)]
\emph{
The \( p \times n  \) matrix \( \PPsi \cc{W}_{1}(x) \) is of full row rank, that is its rows are linearly independent as the Euclidean vectors.}
\end{description}
\begin{remark}\label{Rem_full row rank for 1 scale}
In view of Assumption \( (A2) \) below in Section \ref{section: adaptive procedure} it is sufficient to formulate this assumption only for
\( k=1 \), the positive definiteness of other \( \B{k} \)'s follows automatically.
\end{remark}
\begin{remark}\label{Rem_Gram matrix}
The empirical semi-norm of a function \( g(\cdot) \) given by
\( \| g \|_n^2 = n^{-1} \sum_{i=1}^n g^2(\Xi) \), \( \Xi \in \cc X \), is generated by the ``empirical'' scalar product associating the scalar product in
\( \RRn \): \( \; \langle g,  f \rangle_n = n^{-1} \sum_{i=1}^n g(\Xi) f(\Xi) \; \) with the functions \( g(\cdot) \) and \( f(\cdot) \). Given a weight function
\( s(\cdot)>0 \) one can define in a similar way a weighted empirical scalar product
\begin{equation*}
    \langle g, f \rangle_{n,s} = n^{-1} \sum_{i=1}^n g(\Xi) f(\Xi) s(\Xi)
\end{equation*}
and the corresponding weighted empirical semi-norm. Thus we see that given \( \s(\cdot)>0 \) and a collection of functions
 \( \{ \w{k}{\cdot}(x) \}_{k=1}^K \), the matrices \( n^{-1}\B{k}(x) \) are the Gram matrices of the localized basis functions
 \( \psi_1, \ldots , \psi_p \) centered at \( x \), that is, for any \( k=1, \ldots , K \) we have
\begin{equation*}
    n^{-1}\B{k}(x) = ( \;  \langle \psi_{\nu} (\cdot-x) \sqrt{\w{k}{\cdot}} \, , \;  \psi_{\eta} (\cdot-x) \sqrt{\w{k}{\cdot}} \rangle_{n,\s} \;  )_{1 \le \nu \le \eta \le p },
\end{equation*}
where
\begin{equation*}
     \langle \;  g(\cdot - x)\sqrt{\w{k}{\cdot}} \,  , \;  f(\cdot - x)\sqrt{\w{k}{\cdot}} \; \rangle_{n,\s} = n^{-1} \sum_{i=1}^n g(\Xi- x) f(\Xi- x)\w{k}{i}(x) \s_i^{-2}
\end{equation*}
with \(\s_i = \s(\Xi) >0 \) and \( \w{k}{i}(x) = \w{k}{\Xi}(x) \).
It is well known that any Gram matrix is non-negative definite. Correspondingly, \( \B{k}\succ 0 \) if and only if the rows of \( \PPsi \cc{W}_{1}(x)^{1/2} \) are linearly independent. In view of \( (A2) \) it is sufficient to formulate this assumption only for \( k=1 \). We require slightly more: the rows of \( \PPsi \cc{W}_{1}(x) \)
to be independent. This guarantees that all the variances \( \Var[\mmle{k}] \) are non-degenerated. Indeed, from \eqref{Vk def} below we have  \( \Var[\mmle{k}]=\B{k}^{-1}\tilde{\B{k}} \B{k}^{-1}  \), where \( \tilde{\B{k}}(x) = \PPsi\W{k} \Sigma_0 \W{k}  \PPsi^{\T} \) is a Gram matrix of the same as \( \B{k}(x) \) type, but generated by the scalar products
\( \langle g(\cdot - x)\w{k}{\cdot} \;, f(\cdot - x) \w{k}{\cdot} \rangle_{n,\s}  \)
\end{remark}

\par The formulas in \eqref{MLE_formula} give a sequence of estimators \( \{ \mmle{k}(x) \}_{k=1}^K \). It was noticed in~\cite{Akaike} that in the case when the true data distribution is unknown the QMLE is a natural estimator for the parameter maximizing the expected log-likelihood. That is, for every \( k=1, \ldots, K \), the estimator \( \mmle{k}(x) \) can be considered as an estimator of
\begin{eqnarray}\label{prameter of BPF}
    \bbpf{k}(x) &\eqdef&
     \argmax_{\tta \in \Theta}  \EE \LL \left( \W{k}, \tta \right) \\
     &=& \argmin_{\tta \in \Theta}  (\ff-\PPsi^{\T}\tta)^{\T} \W{k} (\ff-\PPsi^{\T}\tta)\nn
      &=&   \B{k}^{-1} \PPsi \W{k}\ff
        =\B{k}^{-1} \sum_{i=1}^{n}  \Psii f(\Xi) \frac{\w{k}{i}}{ \sigma_i^2}.
\end{eqnarray}
Recall that we do not assume \( \ff = \PPsi^{\T} \tta\) even locally.
It is known from \cite{White} that in the presence of a model misspecification for every \( k \)
the QMLE \( \mmle{k} \) is a strongly consistent estimator for \( \bbpf{k}(x) \), which also is the minimizer of the weighted Kullback-Leibler \cite{KullbackLeibler} information criterion:
\begin{eqnarray*}
  \bbpf{k}(x) &=&
     \argmin_{\tta \in \Theta} \sum_{i=1}^n
        \KL \left( \norm{\ffi}{\s_i} , \norm{\Psii^{\T} \tta}{\s_i}\right) \w{k}{i}(x) \\
  &=& \argmin_{\tta \in \Theta} \sum_{i=1}^n |\ffi - \Psii^{\T} \tta|^2 \frac{\w{k}{i}(x)}{\s_i^2}
\end{eqnarray*}
with \( \KL(P,P_{\tta}) \eqdef \EE_P \big[\log \big(\dd P/\dd P_{\tta}  \big)\big] \). For properties of the Kullback-Leibler divergence see, for example, \cite{Tsybakov}.

\par It follows from the above definition of \( \bbpf{k}(x) \)  and from \eqref{MLE}
that the QMLE \( \mmle{k} \) admits a decomposition into
deterministic and stochastic parts:
\begin{eqnarray}
  & & \mmle{k}
  = \B{k}^{-1}\PPsi \W{k}  (\ff + \Sigma_0^{1/2} \eeps)
  = \bbpf{k} + \B{k}^{-1}\PPsi \W{k}
    \Sigma_0^{1/2} \eeps \label{linearity if quasiMLE}\\
  & & \EE \mmle{k}
    = \bbpf{k},
  \end{eqnarray}
  where \( \eeps \sim \norm{0}{I_n} \). Notice that if \( \ff \equiv \PPsi^{\T}\tta \), then \( \bbpf{k} \equiv \tta \) for any \( k \), and the classical parametric set-up takes place.

\subsection{Adaptive bandwidth selection}
\label{section: adaptive procedure}
Let a point \( x \in \cc X \subset \RRn \), basis \( \{ \psi_j \} \)
and method of localization \( w_{(\cdot , \cdot )  } (x) \) be fixed.
The crucial assumption for the procedure under consideration to work is that the localizing schemes (scales)
\( \cc{W}_{k}(x) = \diag(\w{k}{1}(x) , \ldots , \w{k}{n}(x) ) \) are nested, see Remark~\ref{Rem_Heuristic_explanation_of_the_procedure}. We say that the localizing schemes are nested if for the corresponding matrices the following {\it ordering condition} is fulfilled:
\begin{description}
\item[\( \bb{(A2)} \)]
    \emph{ For any fixed \( x \) and the method of localization with \( w_{(\cdot, \cdot) } (x) \ge 0 \) the following relation holds:
    \begin{equation*}
    \cc{W}_{1}(x) \le  \ldots \le  \cc{W}_{k}(x) \le  \ldots \le  \cc{W}_{K}(x).
            \end{equation*}    }
\end{description}
The inequalities are understood componentwise: for \( 1 \le l\le k \le K \) \( \cc{W}_{l}(x) \le \cc{W}_{k}(x) \) \( \Leftrightarrow \) \( \w{k}{i}(x)- \w{l}{i}(x) \ge 0\) for all \( i= 1, \ldots , n \).
For the kernel smoothing this condition means the following. Given a sequence of bandwidths
\(0 < h_1 <  \ldots < h_k <  \ldots < h_{K} \le 1 \) let \( \w{k}{i}(x) = W( (\Xi-x)h_k^{-1} ) \in [0,1]\) be such that \( W ( u/h_l) \leq W ( u/h_k) \) for any \( 0<h_l < h_k<1 \), and \( W(u) \to 0 \) as \( \| u \| \to \infty \), or even is compactly supported. Also it is intrinsically assumed that, starting from the smallest window, at every step of the procedure every new window contains at least \( p \) new design points.
\par Given the point \( x \in \cc X \),  basis \( \{ \psi_j \} \) and method of localization
\( w_{(\cdot, \cdot)} (x) \), we look for the estimator \(  f_{\adapest} (x) \) of \( f(x) \) having form \eqref{form of the estimator}, where the coefficients \( \adapest^{(j)} (x)  \) are the components
of the estimator
\begin{equation}\label{aadapest}
  \aadapest(x) \eqdef  \mmle{\adapind}(x) =
        (\tilde{\theta}^{(1)}_{\adapind}(x), \ldots, \tilde{\theta}^{(p)}_{\adapind}(x))^{\T},
\end{equation}
corresponding to the adaptive choice of the index \( \adapind \in \{ 1, \ldots, K \} \), i.e. to the choice of the scale. One should keep in mind that \( \adapind \) is a random variable taking values in \( \{ 1, \ldots, K \} \).

\par The selection of \( \aadapest(x) \) from \( \{ \mmle{k}(x)  \} \), \( k = 1, \ldots , K \), can be done by application of the Lepski~\cite{Lep1990} method to comparing of the maximized log-likelihoods
\( \LL(\W{k},\mmle{k}) \). This is the idea of the \emph{fitted local likelihood} (FLL) technique suggested in \cite{KatkSpok}. More precisely, to describe the test statistic, define for any
 \( \tta \), \( \tta' \in \Theta\) the corresponding log-likelihood ratio:
  \begin{equation}\label{log-likelihood ratio}
     \LL(\W{k},\tta, \tta') \eqdef \LL(\W{k},\tta) -\LL(\W{k},\tta'),
  \end{equation}
with \( \LL(\W{k},\tta) \) defined by \eqref{log-likelihood}.

For every \( l = 1, \ldots, K \), the ``fitted'' log-likelihood ratio is defined as follows:
\begin{equation*}
    \LL(\W{l},\mmle{l}, \tta') \eqdef \max_{\tta \in \Theta} \LL(\W{l},\tta, \tta').
  \end{equation*}
By Lemma \ref{Th. Spokoiny_fitted likelihood}, for any scale index \( l \) and parameter vector \( \tta \) this quantity is quadratic in \( \tta \):
\begin{equation*}
    2 \LL(\W{l}, \mmle{l}, \tta) = (\mmle{l} - \tta)^{\T}
        \B{l} (\mmle{l} - \tta).
\end{equation*}
This prompts, see Remark~\ref{Rem_Heuristic_explanation_of_the_procedure}, to use the \emph{FLL-statistics}:
 \begin{eqnarray}\label{Tlk}
    T_{lk} &\eqdef&
                2 \LL(\W{l}, \mmle{l}, \mmle{k})\nn
           & = &
                (\mmle{l} - \mmle{k})^{\T}  \B{l} (\mmle{l} - \mmle{k})\;,\;\;\;\; l < k.
 \end{eqnarray}
In the algorithm \eqref{adaptind} the scale corresponding to \( k=1 \) is assumed to provide \( \B{1} \succ 0 \) and to be sufficiently small assuring nonsignificant deviation of the parametric fit from the true model and \( k=1 \) is always accepted. Then the adaptive index \( \adapind \) is selected by Lepski's selection rule with the FLL test statistics \( \{ T_{lm} \} \):
 \begin{equation}\label{adaptind}
    \adapind = \max \left\{k \leq K :   T_{lm}  \leq \z_l, \, 1 \le l < m \leq k    \right\}.
 \end{equation}
Finally put \( \aadapest = \mmle{\adapind} \).

\par The procedure \eqref{adaptind} involves parameters \( \z_1 , \ldots , \z_{K-1} \). As in the classical Lepski procedure, c.f.~\cite{Lep1990} and \cite{LepMamSpok97}, the inequalities in~\eqref{adaptind} control the risk of estimators for the case of dominating bias. The opposite case of the negligible w.r.t. the noise bias can be easily controlled in view of the Wilks-type result of Proposition~\ref{Wilks theorem}, c.f. Corollary~\ref{param. risk bounds} and Remark~\ref{Rem_Heuristic_explanation_of_the_procedure}:
\begin{equation}\label{bound for the expected fitted
log-likelihood ratio}
    \EE |2 \LL( \W{k}, \mmle{k}, \bbpf{k} )  |^r  \le  C(p,r)
\end{equation}
with the constant \( C(p,r) \) explicitly given by \eqref{C(p,r)} in Appendix.
\par Let \( \aadapest_{k} \) denote the last accepted estimate after the first \( k \) steps of the procedure:
  \begin{equation}\label{the last accepted}
     \aadapest_{k} \eqdef \mmle{\min \{k,\adapind\}}.
  \end{equation}
 Suppose at this step that the critical values \( \z_1 , \ldots , \z_{K-1} \) have being fixed satisfying the following set of \( K-1 \) conditions:
 \begin{definition}\textbf{Propagation conditions (PC)}
\par \noindent
Let for a given \( \alpha \in (0,1] \) and \( r>0 \) the critical values \( \z_1, \ldots ,\z_{K-1} \) satisfy
 \begin{equation}\label{PC}
    \EE_{0, \Sigma} |(\mmle{k} - \aadapest_{k})^{\T} \B{k} (\mmle{k} - \aadapest_{k})|^{r}
        \leq \alpha  C(p,r)\; \; \;        \text{for all}\;\; k=2, \ldots, K,
 \end{equation}
where \( C(p,r) \) is defined by \eqref{C(p,r)} and \( \EE_{0, \Sigma} \) stands for the expectation w.r.t.
the measure \( \norm{0}{\Sigma} \).
 \end{definition}
\begin{remark}\label{Rem_pivotality}{\it ``True'' value of \( \tta \).}
Lemma \ref{Pivotality property} from Section \ref{section:Auxiliary results} shows that in the ``no bias'' situation the Gaussian distribution
provides a nice pivotality property: the actual value of the parameter \( \tta \) is not important for
the risk of adaptive estimate, so one can put \( \tta = 0 \) in \eqref{PC}.
\end{remark}
\begin{remark}\label{Rem_practical_choice of CVs}{\it Calculation of the thresholds.}
Clearly at any step \( k \le K \) of the algorithm the ``current value'' of the adaptive estimator \( \aadapest_{k} \) depends on the thresholds \( \z_1, \ldots ,\z_{k-1} \). The theoretical aspects related to the heteroscedasticity of model and incorrectly known variance is the focus of the present paper. Thus we do not detail the practical aspects of the thresholds calibration only mentioning that in practice this can be done by Monte Carlo simulations under the known ``parametric'' measure \( \norm{0}{\Sigma} \). Moreover one needs to calculate them only once. For detailed consideration of the practical aspects of the calibration as well as for the computational results see~\cite{SV} or~\cite{KatkSpok} focused on the image denoising by local constant fitting, where the similar idea was proposed. Demo-versions of the software are available on the web page http://www.cs.tut.fi/\textasciitilde lasip/.
\end{remark}

\begin{remark}\label{Rem_r_alpha}{\it Loss power \( r \) and ``confidence'' level \( \alpha \).}
The choice of the parameters \( \alpha \) and \( r \) is free and depends only on desired accuracy results and procedure performance. The basic oracle result of Theorem~\ref{oracle result} is formulated in terms of polynomial loss function with index \( r/2 \). Therefore the choice of \( r \) in the PC's determines the final risk bounds. The constant \( \alpha \) appears in the second order term of the bound.

\par A detailed explanation of the heuristics behind the PC's and the role of the parameters \( r \) and \( \alpha \) from the hypothesis testing point of view is given in~\cite{SV}, pp.~2789-2790. Below in Remark~\ref{Rem_Heuristic_explanation_of_the_procedure} we present other heuristics for the procedure and PC's, also explaining why \( \alpha \le 1 \). Here we just mention that the result of Proposition~\ref{Prop_theoret.choice_CV} shows that up to the constants the critical values \( \z_1 , \ldots , \z_{K-1} \) are of the form \( \z_k = C_1 r (K-k) + C_2 \log(K/\alpha) + C_3 \). Therefore the high value of \( r \) along with small \( \alpha \) enlarge \( \z_k \)'s and make the procedure less sensitive to deviations of the parametric fit from the true model resulting in acceptance of a larger smoothing window. Small \( r \) and \( \alpha \) close to one may result in a less stable performance of the procedure and undersmoothing. The free choice of these parameters allows a practical adjustment of the procedure to a particular data set.
\end{remark}

\begin{remark}\label{Rem_Heuristic_explanation_of_the_procedure}{\it Some heuristics behind the procedure.}
Let us give an explanation in the spirit of the example with two H\"older classes (naturally nested w.r.t. the smoothness parameters!) from \cite{Lep1990}, p.2. Let we have only two scales \( \cc W_1(x) \le \cc W_2(x) \) and, correspondingly, two MLE estimators \( \mmle{1}(x) \) and \( \mmle{2}(x) \). The aim is to select automatically from \( \{ \mmle{1}, \mmle{2} \} \). Assume that the noise is known and that either the parametric model \eqref{PA} is true ``globally'', i.e. on \( \cc W_2 \) and consequently (due to (A2)) on \( \cc W_1 \), either \eqref{PA} is satisfied only on \( \cc W_1 \). Two wrong choices are possible:
\begin{description}
  \item[(I)] \( \aadapest = \mmle{1} \) in the global parametric situation when the correct estimator is \( \mmle{2} \);
  \item[(II)] \( \aadapest = \mmle{2}\) when the parametric model is true only on \( \cc W_1 \) and the correct estimator is \( \mmle{1} \).
\end{description}
These two situations are highly asymmetric.

\par Consider (I). Here \( \bbpf{1} = \bbpf{2}=\tta \) and
\( \EE = \EE_{\tta} \), that is \( \EE_{\tta} \mmle{1} = \bbpf{1} =  \bbpf{2}= \EE_{\tta} \mmle{2} \). We have accepted the worst estimator corresponding to the smaller amount of data \( \aadapest = \mmle{1} \) with larger variance. Since \( \cc W_1(x) \le \cc W_2(x) \), by~\eqref{Vk bound} we have \( \Var \mmle{1} = \B{1}^{-1} \succeq \B{2}^{-1} = \Var \mmle{2}\) for the binary weights; for the non-binary weights in \( [0,1] \) \( \Var \mmle{l} \preceq \B{l}^{-1} \), \( l=1,2 \), and the matrices \( \B{l}^{-1} \) serve as monotonized bounds for the variances. Adding and subtracting
\( \LL(\W{2}, \mmle{2}) \) we get
\begin{equation*}
    \LL(\W{2}, \aadapest, \tta)\ind\{ \aadapest = \mmle{1} \}= \LL(\W{2}, \mmle{1}, \mmle{2}) + \LL(\W{2}, \mmle{2}, \tta).
\end{equation*}
Let \( r=1 \). The risk of this log-likelihood ratio is
\begin{equation*}
    \EE_{\tta}|2\LL(\W{2}, \aadapest, \tta)|\ind\{ \aadapest   =   \mmle{1} \}
    \le
    \EE_{\tta}|2\LL(\W{2}, \mmle{1}, \mmle{2})| + \EE_{\tta}|2\LL(\W{2}, \mmle{2}, \tta)|.
\end{equation*}
The second term of the RHS is bounded with \( C(p,1) \) by Corollary~\ref{param. risk bounds}. On the contrary, the first term related to the ``pure noise'' (the value \( \tta \) cancels in \( \mmle{2} - \mmle{1} \)) by the second statement of Lemma~\ref{param. risk bounds polinom} can be much larger than \( C(p,1) \). However, because the distribution of this quantity does not depend on the unknown parameter \( \tta \), its risk can be easily controlled by the choice of the threshold \( \z_1 \). Thus we have arrived at the PC: \( \z_1 \) should provide \( \EE_{\tta}|2\LL(\W{2}, \mmle{1}, \mmle{2})| \le \alpha C(p,1) \) with some \( \alpha \in (0,1] \).

\par In a general case at this place one needs exponential inequalities to bound the large deviations of the stochastic term in the ``no noise'' situation. For analysis of large deviations of a contrast function related to the considered here approach see \cite{Golubev_Spokoiny}.

\par Turn now to (II). Here \( \bbpf{1} \neq \bbpf{2} \). Similarly to the previous case we have
\begin{equation*}
    \EE|2\LL(\W{1}, \aadapest, \bbpf{1})|\ind\{ \aadapest   =   \mmle{2} \}
    \le
    \EE|2\LL(\W{1}, \mmle{2}, \mmle{1})| + \EE|2\LL(\W{1}, \mmle{1}, \bbpf{1})|
\end{equation*}
and as in (I) the second term of the RHS is bounded with \( C(p,1) \). But one can say nothing about the first term and the only way to control it is the procedure: we say that the choice \( \aadapest   =   \mmle{2} \) is acceptable in this situation if \( 2\LL(\W{1}, \mmle{1}, \mmle{2}) \le \z_1 \), where \( \z_1  \) is the threshold fixed by the PC. To choose from more than two estimators the selection rule at every step \( k \) accepts the estimator \( \mmle{k} \) if and only if  \( 2\LL(\W{l}, \mmle{l}, \mmle{k}) \le \z_l  \) for all \( l < k \) with the proviso that \( \mmle{k-1} \) had been accepted at the previous step of the procedure.

\par Note also that exactly this part of the procedure can cause the well-known oversmoothing effect of the Lepski-type procedures, because one admits oversmoothing in the range of threshold. The threshold corresponding to the oracle scale presents also in the leading term of the risk, c.f. Theorem~\ref{oracle result}. That is why it is so important to select the smallest possible sequence of thresholds and it is shown in \cite{SV} p. 2791 that the PC's provide such a sequence. However, to fix the thresholds by simulations as in \cite{SV} the exact knowledge of the noise is required. This explains the interest of the author to the noise misspecification and generalization of the propagation approach to this set-up.
\end{remark}

\section{Theoretical study}
\par In order to infer on the admissible level of misspecification for ``model'' covariance matrix from \eqref{PA} we need to introduce a parameter \( \delta \) reflecting the relative variability in errors:
\begin{description}
\item[\( \bb{(A3)} \)]
\emph{
There exists \( \delta \in [0,1) \) such that
\begin{equation*}
1-\delta \le \sigma_{0,i}^2/\sigma_{i}^2 \le 1+ \delta \; \;\; \text{for all} \;\; \; i = 1, \ldots, n.
\end{equation*}
}
\end{description}
\begin{remark}
Clearly, the value of \( \delta \) is not available. This parameter is used to trace the influence of the erroneously known noise. The procedure given by \eqref{Tlk}, \eqref{adaptind} and \eqref{PC} does not require knowledge of \( \delta \) or of the true covariance matrix \( \Sigma_0 \).
\end{remark}

\subsection{Upper bound for the critical values}
\par  For any real symmetric matrices \( A \) and \( B \) we write \( A \preceq B \) if \( \vartheta^{\T} A \, \vartheta  \leq \vartheta^{\T} B \, \vartheta  \) for all vectors \( \vartheta \), or, equivalently, if and only if the matrix \(  B -  A \) is nonnegative definite. Assuming \( {(A3)} \), the true covariance matrix \( \Sigma_0 \preceq \Sigma (1+\delta) \), and the variance of the estimate \( \mmle{k} \) is bounded with \( \B{k}^{-1} \):
\begin{eqnarray}
     V_{k} \eqdef \Var \mmle{k}
        &=&    \B{k}^{-1}\PPsi\W{k} \Sigma_0 \W{k}  \PPsi^{\T}\B{k}^{-1} \label{Vk def}
\\        &\preceq &  (1+\delta)\B{k}^{-1}\PPsi\W{k} \Sigma  \W{k} \PPsi^{\T}\B{k}^{-1}\nn
        &=&  (1+\delta) \B{k}^{-1}\PPsi  \Sigma^{-1/2} \cc W_{k}^2 \Sigma^{-1/2}  \PPsi^{\T}\B{k}^{-1}\nn
        &\preceq & (1+\delta) \B{k}^{-1}\PPsi  \Sigma^{-1/2} \cc W_{k} \Sigma^{-1/2} \PPsi^{\T}\B{k}^{-1}\nn
       &=& (1+\delta) \B{k}^{-1} \PPsi \W{k} \PPsi^{\T} \B{k}^{-1} \nn
       &=& (1+\delta) \B{k}^{-1}.\label{Vk bound}
       \end{eqnarray}
The last inequality follows from the observation that all entries of the diagonal ``weight''
matrix \( \cc{W}_{k} \) do not exceed one, implying \(\cc W_{k}^2 \preceq \cc W_{k} \). The strict equality
 takes place if \( \{ \w{k}{i} \}  \in \{ 0,1 \}\) and the noise is known, i.e. if \( \delta = 0 \).
To justify the procedure it is necessary to show that the critical values fixed by \((PC)\) are finite.
This will be obtained under the following assumption:
\begin{description}
\item[\( \bb{ (A4) } \)]
\emph{
Let for some constants \( u_0 \) and \( u \) such that \( 1< u_0 \le  u \)
for any \( 2\le k\le K \) the matrices \( \B{k} \) satisfy
\begin{equation*}
   u_0 I_{\p} \preceq  \B{k-1}^{-1/2} \, \B{k} \, \B{k-1}^{-1/2} \preceq u I_{\p}
\end{equation*}
}
\end{description}
 \begin{remark}\label{Rem_connection_Bk_with_bandwidth}
 In the ``one dimensional case'' \( \p=1 \), that is for the local constant fitting, the ``matrix''
 \( \B{k} = \sum_{i=1}^{n} \w{k}{i} \sigma_i^{-2} \ge  \B{k-1} \) is just a weighted ``local sample size''. Assume for simplicity that \( \sigma_i^{2} \equiv \sigma^2 \), the weights are rectangular kernels \( \w{k}{i}(x) = \ind\{ |\Xi - x| \le h_k/2 \} \) and the design is equidistant. Then for \( n \) sufficiently large
 \begin{equation*}
     \frac{1}{n} \B{k} = \frac{1}{n\sigma^2} \sum_{i=1}^{n} \ind\{ |\frac{i}{n} - x| \le \frac{h_k}{2} \} \approx  \frac{h_k}{\sigma^2},
 \end{equation*}
and Assumption \( (A4) \) with \( u_0 = u \) means that the bandwidths grow geometrically: \( h_k =u h_{k-1} \).
 \end{remark}
Now we are able to demonstrate the finiteness of the critical values.
\begin{proposition}\label{Prop_theoret.choice_CV} {\it \textbf{Theoretical choice of the critical values}}
\label{upper bound}
\par Assume \( (A1) - (A2) \) and \( (A4) \).
The adaptive procedure defined by \eqref{Tlk}, \eqref{adaptind} and \eqref{PC} is well defined in the sense that the choice of the critical values of the form
\begin{equation}
    \z_k = \frac{4}{\mu} \left\{ r(K-k) \log u + \log{(K/\alpha)} - \frac{p}{4} \log(1-4\mu )
            - \log(1-u^{-r}) + \bar{C}(p,r) \right\} \label{zk}
\end{equation}
provides the conditions \eqref{PC} for all \( k\le K \).
Particularly,
\begin{equation}
    \EE_{0, \Sigma} |(\mmle{K} - \aadapest)^{\T}
        \B{K} (\mmle{K} - \aadapest)|^{r}
        \leq \alpha  C(p,r).
 \end{equation}
In \eqref{zk} \( \mu \in (0,1/4) \) is an arbitrary constant, \( u>1 \) is given by Assumption \( (A4) \), \( r>0 \) and \( \alpha \in (0,1] \) are from the PC's, and
 \begin{equation*}
     \bar{C}(p,r) =
     \log\left\{ \frac{2^{2r} [\Gamma(2r + p/2) \Gamma(p/2)]^{1/2}}{\Gamma(r + p/2)} \right \}.
 \end{equation*}
\end{proposition}
The proof is given in Section \ref{proof of CV}.
\begin{remark}\label{Rem_dependence_of_th_on_param}
\par {\it Dependence of the thresholds on the parameters} \( r \) {\it and} \( \alpha \) in connection with the performance of the procedure is discussed in Remark~\ref{Rem_r_alpha}.
\par \emph{Dependance on the number of scales.} For kernel estimators, c.f. Remark~\ref{Rem_connection_Bk_with_bandwidth}, Assumption \( (A4) \) essentially means a geometrical grid of bandwidths \( h_{k-1} = u^{-1} h_k \) implying \( h_1 = u^{-(K-1)} h_K \). Thus \( (K-1) \log u = \log(h_K/h_1) \), where \( \log u \) is a fixed constant, say equal to \( \log 2 \). Since \( h_K \le 1 \) and \( h_1 \ge 1/n \), the number of scales is at most of order \( \log n \), that is \( K \asymp \log(h_K/h_1) \le \log n \) and is related to the ``adaptive factor'' to pay for the pointwise adaptation, c.f. (2.11) in \cite{LepSpok97} p.~2518 and the discussion therein. The leading term in~\eqref{zk} is \( const.(K-k) \) and it shows that the thresholds \( \z_k \) linearly decrease in \( k \) providing stability of the procedure at the first steps and sensitivity to deviations of the parametric fit from the true model at the further steps of the algorithm. The thresholds are at most of order \( \log n \) and this ``log'' disappears at the ``last point'' \( k=K \). That is if the parametric assumption is true, there is no ``log-payment'', c.f. Remark~\ref{Rem_comments on the basic oracle result}.
\end{remark}

\subsection{Quality of estimation in the nearly parametric case}
\label{nearly parametric case}
The critical values of the procedure \( \z_1, \ldots , \z_{K-1} \) were selected by the propagation conditions~\eqref{PC} under the measure \( \norm{\tta}{\Sigma} \)that is probably not confirmed by the data. Let now the maximizers of the expected local log-likelihoods \( \bbpf{1}, \ldots, \bbpf{k}  \) are only approximately equal, say to \( \tta \), up to some \( k\leq K \) and the covariance matrix is \( \Sigma_0 \). The meaning of ``approximately equal'' will be explained below.

\par The aim is to justify the use of the critical values in this situation. For this purposes we study the discrepancy between the joint distributions of linear estimators \( \mmle{1}, \ldots , \mmle{k} \) for \( k = 1, \ldots , K \) under the ``no bias'' assumption corresponding to the distributions with mean \( \bbpf{1} = \cdots = \bbpf{k} = \tta\) and possibly incorrectly specified covariance matrix \( \Sigma \), and in the general situation with \( \bbpf{1} \ne \cdots \ne \bbpf{k} \) and covariance \( \Sigma_0 \). Denote the expectations w.r.t. these measures by \(\EE_{\tta, \Sigma} := \EE_{k,\tta, \Sigma} \) and \( \EE_{\ff, \Sigma_0} := \EE_{k,\ff, \Sigma_0} \) respectively and the \( p\times k \) matrix of the first \( k \) estimators and the expectations correspondingly by
\begin{eqnarray*}
  \tilde\TTa_k &\eqdef& (\mmle{1}, \ldots, \mmle{k}), \\
  \TTa^*_k &\eqdef& \EE_{\ff, \Sigma_0} \tilde \TTa_k =(\bbpf{1}, \ldots, \bbpf{k}), \\
  \TTa_k &\eqdef& \EE_{\tta, \Sigma} \tilde \TTa_k= (\tta, \ldots, \tta).
\end{eqnarray*}
Let \( A \otimes B \) stand for the Kronecker product of matrices \( A = (a_{i,j})_{1 \le i \le m , 1 \le j \le n} \) and \( B \) defined as
\begin{equation*}
    A \otimes B =\left(
                   \begin{array}{cccc}
                     a_{11}B & a_{12}B & \cdots & a_{1n}B \\
                     a_{21}B & a_{22}B & \cdots & a_{2n}B \\
                     \cdot   & \cdot   & \cdots & \cdot \\
                     a_{m1}B & a_{m2}B & \cdots & a_{mn}B \\
                   \end{array}
                 \right).
\end{equation*}

Denote the \( pk \times pk \) covariance matrices of \( \vec \tilde \TTa_k^{\T} = (\mmle{1}^{\T}, \ldots , \mmle{k}^{\T}) \in \RR^{pk}\) by
\begin{eqnarray}
 \SSigma_k    &\eqdef&  \Var_{\tta, \Sigma}[\vec \tilde \TTa_k] =\DD_k (J_k \otimes \Sigma) \DD_k^{\T}, \label{def SSigma}\\
\SSigma_{k,0}   &\eqdef&  \Var_{\ff, \Sigma_0} [\vec \tilde \TTa_k] =\DD_k (J_k \otimes \Sigma_0) \DD_k^{\T} \label{def SSigma0},
\end{eqnarray}
where the matrix \(J_k \) is a \( k \times k \) matrix with all its elements equal to \(1\), and the
 \( pk \times nk \) block diagonal matrix \( \DD_k \) is defined as follows:
\begin{eqnarray}
 \DD_k    &\eqdef&    D_1 \oplus \cdots \oplus D_k = \diag(D_1, \ldots, D_k), \label{def DD_K}\nn
  D_l     &\eqdef&    \B{l}^{-1} \PPsi \W{l}, \;\;\; l = 1, \ldots ,k.
\end{eqnarray}
By Lemma \ref{simidefinitness of SSigma } from Section \ref{section:Auxiliary results} under Assumption \( (A3)\) with the same \( \delta \) the similar relation holds for the covariance matrices \( \SSigma_{k}   \) and \( \SSigma_{k,0} \) of the sets of linear estimators:
\begin{equation}\label{multi cond sigma}
  (1-\delta) \SSigma_{k} \preceq  \SSigma_{k,0} \preceq (1+\delta)  \SSigma_{k}\;,\;\; k \le K.
\end{equation}
\par In spite of by Lemma \ref{MGF for joint distribution} the moment generating function of \( \vec \tilde \TTa_K \) has the form corresponding to the multivariate normal distribution this representation makes sense only if \( \SSigma_K \) is nonsingular. Notice that \( \rank(J_K \otimes \Sigma) =n \). From \( J_K \otimes \Sigma \succeq 0\) it follows only that \( \SSigma_{K} \succeq 0\), similarly, \(  \SSigma_{K,0} \succeq 0 \). However, without any additional assumptions it is easy to show, see Lemma \ref{nonsingularity of SSigma rectang}, that for rectangular kernels \( \SSigma_{K} \succ 0 \). On the other hand, due to \eqref{multi cond sigma}, it is enough to require nonsingularity only for the matrix \( \SSigma_K \) corresponding to the approximate model \eqref{PA}, and its choice belongs to a statistician. In what follows we assume that \( \SSigma_{K} \succ 0\).

\par Denote by \( \P_{\tta, \Sigma}^k = \norm{\vec\TTa_k}{\SSigma_k} \) and by \( \P_{\ff, \Sigma_0}^k = \norm{\vec\TTa^*_k}{\SSigma_{k,0}} \), \( k=1, \ldots, K \), the distributions of \( \vec \tilde \TTa_k \) under the assumption that the parametric model \eqref{PA} is true up to the scale \( k \) and under the assumption that nonparametric model \eqref{true model} takes place. Denote also the Radon-Nikodym derivative by
\begin{equation}\label{def Zk}
    Z_k \eqdef \frac{\dd \P_{\ff, \Sigma_0}^k}{\dd \P_{\tta, \Sigma}^k}.
\end{equation}
Then Lemma \ref{KL for joint distr} gives the Kullback-Leibler divergence between these measures:
\begin{eqnarray}
  & &2\KL(\P_{\ff, \Sigma_0}^k,\P_{\tta, \Sigma}^k) \eqdef 2\EE_{\ff, \Sigma_0}
  \log (Z_k) \label{KL joint}\nn
  &=&  \Delta(k)  +\log\bigg(\frac{\det \SSigma_k}{\det \SSigma_{k,0}} \bigg)
            + \tr (\SSigma_k^{-1} \SSigma_{k,0}) -pk,
\end{eqnarray}
where
\begin{eqnarray}
 \Delta(k) & \eqdef & b(k)^{\T} \SSigma_k^{-1} b(k)  \label{def_Delta(k)} \\
  b(k) &\eqdef& \vec \TTa^*_k - \vec \TTa_k \label{def_b(k)} .
\end{eqnarray}
If there would be no any ``noise misspecification'', i.e. if \( \delta \equiv 0 \) implying
\( \Sigma = \Sigma_0 \), then
\( \Delta(k) = b(k)^{\T} \SSigma_k^{-1} b(k) = 2 \KL(\P_{\ff, \Sigma}^k,\P_{\tta, \Sigma}^k)\). Under Assumption \( (A2) \), the quantity \( \Delta(k) \) grows with \( k \), so following the terminology suggested in \cite{SV}, we
introduce the \emph{small modeling bias condition}:
\begin{description}
\item[\( \bb{(SMB)} \)]
\emph{
Let for some \( k \le K \) and \( \tta \) exist a finite constant \( \Delta  \ge 0 \) such that \( \Delta(k) \le  \Delta. \)
}
\end{description}
Monotonicity of \( \Delta(k) \) and \( (SMB) \) immediately imply that
\begin{equation*}
   \sup_{1 \le l \le k} \Delta(l) \le \Delta.
\end{equation*}
Relation \eqref{multi cond sigma} yields \( - p k \delta \le  \tr (\SSigma_k^{-1} \SSigma_{k,0}) -p k \le  p k \delta \). Thus the statement of Lemma~\ref{KL for joint distr} gives a bound for the Kullback-Leibler divergence in terms of~\( \delta \):
\begin{eqnarray}
  -\frac{pk}{2} \log(1+ \delta) + \frac{\Delta(k)}{2}  -\frac{ pk \delta}{2} &\le& \KL(\P_{\ff, \Sigma_0}^k,\P_{\tta, \Sigma}^k) \label{bounds for KL} \\ \nonumber
   &\le& -\frac{pk}{2} \log(1- \delta) + \frac{\Delta(k)}{2}  +\frac{ pk \delta}{2}.
\end{eqnarray}
Moreover, if \( \delta = \delta(n) \) and  \( \delta(n) \to 0+\) as \( n \to \infty \)
\begin{equation}\label{bounds for KL asymptotics}
     \Delta(k) -2 pk\delta + o(\delta)
    \le
    2\KL(\P_{\ff, \Sigma_0}^k,\P_{\tta, \Sigma}^k)
    \le
    \Delta(k) + 2 pk\delta + o(\delta).
\end{equation}
This means that, if for some \( k \) Assumption \( (SMB) \) is fulfilled and \( \delta = O(1/K )  \), then the Kullback-Leibler divergence between the measures\( \P_{\tta, \Sigma}^k  \) and  \( \P_{\ff, \Sigma_0}^k  \) is bounded by a small constant.
\par Now one can state the crucial property for obtaining the final oracle result.
\begin{theorem}\label{Propagation result theorem} \textbf{Propagation property}
   \par Assume \( (A1) - (A4) \) and \( ({PC}) \). Then for any \( k \le K \) the following upper bounds hold:
    \begin{eqnarray*}
    && \EE|(\mmle{k} - \tta)^{\T} \B{k} (\mmle{k} - \tta)|^{r/2}\\
       & \le &  C(p,r)^{1/2} (1+\delta)^{pk/4}(1-\delta)^{-3pk/4}
        \exp\left\{ \varphi(\delta)
        \frac{\Delta(k)}{2(1-\delta)}\right\},\\
       && \EE|(\mmle{k} - \aadapest_{k})^{\T} \B{k} (\mmle{k} - \aadapest_{k})|^{r/2}\\
       & \le & (\alpha C(p,r))^{1/2} (1+\delta)^{pk/4}(1-\delta)^{-3pk/4}
        \exp\left\{ \varphi(\delta)
        \frac{\Delta(k)}{2(1-\delta)}\right\},
 \end{eqnarray*}
 where \(\varphi(\delta) \eqdef
  \begin{cases}
1 & \mathrm{for \; homogeneous \;errors,} \\
 \frac{2(1+\delta)}{(1-\delta)^2} -1  &\mathrm{otherwise}.
\end{cases} \)
\par Here \( \mmle{k} = \mmle{k}(x) \) is the QMLE defined by \eqref{MLE}, \( \tta \) is the parameter from \eqref{PA}, \( \aadapest_{k}(x) = \mmle{\min{k, \adapind}}(x) \) is the adaptive estimate at the \( k \)th step of the procedure, \( C(p,r) \) is the constant from the PC's defined in \eqref{C(p,r)} and \( p \) is the number of basis functions used for the linear fitting.
\end{theorem}
The proof is given in Subsection \ref{proof propagation}.
\begin{remark}\label{Rem_relative_noise_error}
Bounds \eqref{bound for Zk homog} and \eqref{bound for Zk} obtained in the proof of the theorem (Section~\ref{proof propagation}) give a condition on the relative error in the noise misspecification. Let \( \delta = \delta(n) \to 0+ \) as \( n \to \infty \). Then for every \( k \le K \)
\begin{equation*}
    \varphi(\delta)
        \frac{\Delta(k)}{1+\delta} -2pk\delta +o(\delta) \le
        \log \EE_{\tta, \Sigma} [Z_k^2]\le
       \varphi(\delta)
        \frac{\Delta(k)}{1-\delta} +2pk\delta + o(\delta)
\end{equation*}
with \( Z_k \) defined by \eqref{def Zk}. This bound implies, up to the additive constant
\( 0.5 \log\big( \alpha C(p,r) \big)\), the same asymptotic behavior for the logarithm of the risk of adaptive estimate \( \log \EE \| \B{k}^{1/2} (\mmle{k} - \aadapest_{k}) \|^r \) at each step of the procedure. Because by \( (SMB) \) the quantity \( \Delta(k) \) is supposed to be bounded by a small constant, and \( K \) is of order \( \log n \), see Remark~\ref{Rem_dependence_of_th_on_param}, the expectation \( \EE_{\tta, \Sigma} [Z_k^2] \) is small if \( \delta = O ( 1/\log n) \) and, consequently, the risk \( \EE \| \B{k}^{1/2} (\mmle{k} - \aadapest_{k}) \|^r \) is bounded, c.f. \eqref{PC}. This means that for a plug-in estimator of the variance only the logarithmic in sample size quality is needed. This observation is of particular importance, since it is known from~\cite{Spokoiny variance} that over classes of functions with bounded second derivative the rate \( n^{-1/2} \) of variance estimation is achievable only for the dimension \( d \le 8 \).
\end{remark}
\begin{remark}
The propagation property guarantees that the adaptive procedure does not stop with large probability while \( \Delta(k) \) is small, i.e., under \( (SMB) \), and if the relative error \( \delta \) in the noise is sufficiently small.
\end{remark}

\subsection{Quality of estimation in the nonparametric case: the oracle result}
\label{subsection:oracle result}
Define the {\it oracle index} as the largest index \( k \le K \) such that \( (SMB) \) holds:
\begin{equation}\label{oracle index}
    k^* \eqdef \max \{ k \le K : \Delta(k) \le \Delta \}.
\end{equation}

\begin{theorem}\label{oracle result}
Let \( \Delta(1) \le \Delta \), i.e. the first estimator is always accepted in the testing procedure.
Let \( \SSigma_{K} \succ 0\) and \( k^* \) be the oracle index. Let \( \mmle{k^*} \) be the nonadaptive estimator defined by \eqref{MLE_formula} corresponding to \( k^* \) and \(  \aadapest \) be an output of the procedure \eqref{Tlk} -- \eqref{adaptind}. Then under \( (PC) \) and assumptions \( (A1) - (A4), (SMB) \) and   the risk between the adaptive and oracle estimators is bounded with the following expression:
\begin{eqnarray}\label{oracle bound statement of the Th}
       && \EE|(\mmle{k^*} - \aadapest)^{\T} \B{k^*} (\mmle{k^*} - \aadapest)|^{r/2}\\
       & \le &   \z_{k^*}^{r/2} +
       (\alpha C(p,r))^{1/2} (1+\delta)^{pk^*/4}(1-\delta)^{-3pk^*/4}
        \exp\left\{ \varphi(\delta)
        \frac{\Delta}{2(1-\delta)}\right\} ,\nonumber
\end{eqnarray}
where \( \varphi(\delta) \) is as in Theorem \ref{Propagation result theorem} and \( C(p,r) \) is the constant from the PC's defined in \eqref{C(p,r)}.
\end{theorem}
\begin{remark}\label{Rem_comments on the basic oracle result}
The second term in the RHS of \eqref{oracle bound statement of the Th} is bounded with a constant with the proviso that \( \delta = O(1/\log n) \), see Remark~\ref{Rem_relative_noise_error}, and the leading term is \( \z_{k^*}^{r/2} \) that by Proposition~\ref{Prop_theoret.choice_CV} has the form \( \z_{k^*} = C_1 r (K-k^*) + C_2 \log(K/\alpha) + C_3 \). The leading term \( K \) is at most of order \( \log n \), see Remark~\ref{Rem_dependence_of_th_on_param}, and is the unavoidable payment for the pointwise adaptation, see Theorem~2 on the lower bound in \cite{Lep1990}. This term cancels if \( k^* = K \), that is when the deviation of the parametric fit from the true model is not significant for all observations. This means that the parametric set-up takes place globally and there is no adaptation involved. The canceling of the \( log \) term at the last point of the range of adaptation in the rate is a common feature of this type procedures, sf. \cite{Lep1990}, \cite{LepMamSpok97}, \cite{LepSpok97}.
\par The LHS of the inequality \eqref{oracle bound statement of the Th} is the mathematical expectation of the oracle log-likelihood ratio \( |2 \LL (\W{k^*}, \mmle{k^*}, \aadapest)|^{r/2} \), or the risk of the difference between the adaptive estimator \( \aadapest \) and its nonadaptive counterpart \( \mmle{k^*} \) normalized by the bound for the variance of the oracle estimator. Recall that by \eqref{Vk bound} in the case of binary weights the matrix \( \B{k^*}^{-1} = \Var{\mmle{k^*}} \); generally we have only \( \Var{\mmle{k^*}} \preceq \B{k^*}^{-1}\). Loosely speaking, the result says that the risk of adaptive estimator is of order of the oracle variance multiplied by the logarithmic factor~\( \z_{k^*} \).
\end{remark}
\begin{proof}
By the definition of the adaptive estimate \( \aadapest = \mmle{ \adapind}  \). Because the events \( \{ \adapind \le k^*\} \) and \( \{ \adapind > k^*\} \) are disjunct, one can write
\begin{eqnarray*}
       && \EE|(\mmle{k^*} - \aadapest)^{\T} \B{k^*} (\mmle{k^*} - \aadapest)|^{r/2}\\
&=& \EE|(\mmle{k^*} - \mmle{ \adapind})^{\T} \B{k^*} (\mmle{k^*} - \mmle{ \adapind})|^{r/2} \ind\{ \adapind \le k^* \} \\
&+&
\EE|(\mmle{k^*} - \mmle{ \adapind})^{\T} \B{k^*} (\mmle{k^*} - \mmle{ \adapind})|^{r/2} \ind \{ \adapind > k^*\}.
\end{eqnarray*}
If \(  \adapind \le k^* \) then \( \aadapest_{k^*} \eqdef \mmle{\min \{ k^*, \adapind \} } = \mmle{ \adapind} \). Thus, to bound the first summand, it is enough to apply Theorem~\ref{Propagation result theorem} with \( k = k^* \).
\par To bound the second expectation, i.e. to bound the fluctuations of adaptive estimate \( \aadapest \) at the steps of the procedure for which the \( (SMB) \) condition is not fulfilled anymore, just notice that for \( \adapind > k^* \) the quadratic form coincides with the test statistic \( T_{k^*,\adapind} \)
\begin{eqnarray*}
  && (\mmle{k^*} - \aadapest)^{\T} \B{k^*} (\mmle{k^*} - \aadapest) \\
  &=& (\mmle{k^*} - \mmle{\adapind})^{\T} \B{k^*} (\mmle{k^*} - \mmle{\adapind})
    = T_{k^*,\adapind}.
\end{eqnarray*}
But the index \( \adapind \) was accepted by the procedure, this means that \( T_{l,\adapind} \le \z_l \) for all \( l < \adapind \) and therefore for \( l=k^* \). Thus
\begin{equation*}
    \EE |(\mmle{k^*} - \aadapest)^{\T} \B{k^*} (\mmle{k^*} - \aadapest)|^{r/2} \ind \{ \adapind > k^*\} \le \z_{k^*}^{r/2}.
\end{equation*}
\end{proof}
\subsection{Componentwise oracle risk bounds}
\label{subsection:estimators of the regression function}
\subsection{Componentwise oracle risk bounds}
\label{subsection:estimators of the regression function}
Theorem \ref{oracle result} provides the oracle risk bound for the adaptive estimator
\( \aadapest(x) = \mmle{\adapind} (x)\) of the parameter vector \( \tta(x) \in \RRp \) corresponding to the estimator \( \hat f_{\aadapest}(x)  \) of type \eqref{form of the estimator}. It is interesting to have a look at the oracle quality of estimation of the components \( \theta^{(1)} , \ldots , \theta^{(p)} \) of the vector \( \tta \) having in mind that the choice of polynomial basis leads to the direct estimation of the value of regression function and the derivatives by the coordinates of \( \aadapest \).

\par Denote by \( LP_k(p-1) \) a local polynomial estimator of order \( p-1 \) corresponding to the \( k \)th degree of localization and by \( LP^{ad}(p-1) \) its adaptive counterpart, i.e. \( LP^{ad}(p-1) \eqdef LP_{\adapind}(p-1) \). If the basis is polynomial and the regression function \( f(\cdot) \) is sufficiently smooth in a neighborhood of \( x \), then \( \aadapest(x) \) is the \( LP^{ad}(p-1) \) of the vector \( ( f(x) , f' (x), \ldots , f^{(p-1)} (x))^{\T} \) of the values of the function \( f \) and its derivatives
at the reference point
\( x \in \RRd \).

\par Now we are going to obtain a similar to the previous section oracle result for the components of the vector \( \aadapest(x) \), particularly for \(  \bb e_j^{\T} \aadapest(x) \), \( j = 1, \ldots, p \), where
\( \bb e_j  \) is the \( j \)th canonical basis vector
in \( \RRp \). As a corollary of this general result in the case of polynomial basis we get an oracle risk bound for  \( LP^{ad}(p-1) \) estimator
of the function \( f \) and its derivatives at the point \( x \).
\par \( LP_k(p-1) \) estimator of \( f^{(j-1)}(x) \) is given by
\begin{eqnarray}\label{def of estimators of f and its derivatives}
  \flej{j-1}{k}{x} &=& e_j^{\T}\mmle{k}(x), \; j= 1, \ldots,p, \\ \nonumber
  \fle{k}{x} &=& \flej{0}{k}{x} = e_1^{\T}\mmle{k}(x).
\end{eqnarray}
Then the adaptive local polynomial estimators are defined as follows:
\begin{eqnarray}\label{def of adaptive estimators of f and its derivatives}
  \adaplpest^{(j-1)}(x) &=& e_j^{\T}\aadapest(x), \; j= 1, \ldots,p, \\ \nonumber
  \adaplpest(x) &=&  e_1^{\T}\aadapest(x).
\end{eqnarray}
Similarly, the adaptive estimators of the function \( f \) and its derivatives corresponding to the \( k \)th step of the procedure are given by
\begin{equation}\label{def of k-th adaptive estimators of f and its deriv}
    \adaplpest_k^{(j-1)}(x) \eqdef e_j^{\T}\aadapest_k(x), \; j= 1, \ldots,p.
\end{equation}
Thus, if the basis is polynomial, the estimator
\( \adaplpest(x) \eqdef \adaplpest^{(0)}(x) \) is the \( LP^{ad}(p-1) \)
estimator of the value \( f(x) \), and \( \adaplpest^{(j-1)}(x) \) with \( j = 2, \ldots, p \)
are, correspondingly, the \( LP^{ad}(p-1) \) estimators of the values of its derivatives.
However the results of Theorems \ref{oracle result} and \ref{oracle result componentwise} hold for any basis satisfying the conditions of the theorems. For the study below we need the following assumptions:
\begin{description}
\item[\( \bb{(A5)} \)]
\emph{There exists a positive finite number  \(  \s_{max}(k)  \) such that for \( i: \Xi \in U_{h_k}(x) \), with the neighborhood of the estimation point \( U_{h_k}(x) \) given by \( \W{k} \) the variances of errors from the parametric (known) model \eqref{PA} are locally uniformly bounded:}
\begin{equation*}
   \s^2_{i} \le \s^2_{max}(k).
\end{equation*}

\item[\( \bb{(A6)} \)]
\emph{Let assumption \( (A5) \) be satisfied. There exists a number \( \Lambda_0>0  \) such that for any \( k=1, \ldots, K \) the smallest eigenvalue \( \lambda_p(\B{k}) \ge  nh_k^d\Lambda_0 \s^{-2}_{max}(k)\) for \( n \) sufficiently large.}
\end{description}
\begin{remark}\label{Rem_Bound for min eig of Bk}
The first assumption is not restrictive at all, since it is about the known variance from the model we use for the construction of estimators. The last assumption is stronger than the requirement \( \B{k} (x) \succ 0   \). Lemmas 1.5, 1.4 in \cite{Tsybakov} shows that this assumption holds for non-negative kernels, which are bounded from below on a set of positive Lebesgue measure. The constant \( \Lambda_0 \) is related to the smallest eigenvalue of the matrix \( \bb B \) from Lemma~\ref{Lm_converg to B}.
\end{remark}
Thus for any \( k= 1, \ldots K \) and for any \( \gamma \in \RRp \) we have
\begin{equation}\label{bound for Bk via the smallest eigenvalue in Rd}
    \gamma^{\T} \B{k}^{-1} \gamma \le \frac{\s^2_{max}(k) }{nh_k^d\Lambda_0} \| \gamma \|^2
                                \le \frac{\bar\s^2_{max}(k) }{nh_k^d\Lambda_0} \| \gamma \|^2,
\end{equation}
where \( \bar\s^2_{max}(k) \eqdef \max_{1 \le l \le k} \s^2_{max}(l) \). Thus we have the following bound:
\begin{lemma}\label{bound for the componentwise differences}
Let \( (A5) \) and \( (A6) \) be satisfied. Then for any \( j = 1, \ldots, p \) and \( k, \, k'= 1, \ldots K \) the following bound holds:
\begin{equation*}
    \left( \frac{nh_k^d\Lambda_0}{\bar\s^2_{max}(k) } \right)^{1/2}
        |\bb e_j^{\T}\mmle{k} - \bb e_j^{\T}\mmle{k'}|
            \le
        \| \B{k}^{1/2} (\mmle{k} - \mmle{k'}) \|.
\end{equation*}
\end{lemma}
\begin{proof}
By \eqref{bound for Bk via the smallest eigenvalue in Rd} taking \( \gamma = \B{k}^{1/2} (\mmle{k} - \mmle{k'})\) we have
\begin{eqnarray*}
  |\bb e_j^{\T}\mmle{k} - \bb e_j^{\T}\mmle{k'}|^2 &\le & \| \mmle{k} -\mmle{k'} \|^2 \\
    &=&
        \|\B{k}^{-1/2} \B{k}^{1/2} (\mmle{k} - \mmle{k'}) \|^2 \\
    &\le &
        \frac{\bar\s^2_{max}(k) }{nh_k^d\Lambda_0}  \| \B{k}^{1/2} (\mmle{k} - \mmle{k'}) \|^2.
\end{eqnarray*}
\end{proof}

\par To obtain the ``componentwise'' oracle risk bounds we need to recheck the ``propagation property''. Firstly, notice that the ``propagation conditions'' \eqref{PC} on the choice the critical values \( \z_1, \ldots, \z_{K-1} \) imply the similar bounds for the components \( \bb e_j^{\T} \aadapest_k(x) \). Recall that \( \aadapest_k \eqdef \mmle{\min\{ k, \adapind \}} \). By \eqref{PC}, Lemma~\ref{bound for the componentwise differences} and the pivotality property from Lemma~\ref{Pivotality property} we have the following simple observation that serves as a componentwise counterpart of PC:
\begin{lemma}\label{PC componentwise} Under the propagation conditions \( (PC) \) for any \( \tta \in \RRp \) and all \( k=2, \ldots, K \) we have:
    \begin{eqnarray*}
     \left( \frac{nh_{k}^d\Lambda_0}{\bar\s^2_{max}(k) } \right)^{r}
               \EE_{\tta, \Sigma} |\bb e_j^{\T}\mmle{k}(x) - \bb e_j^{\T}\aadapest_k(x)|^{2r}
  &\le&  \EE_{0, \Sigma} \| \B{k}^{1/2} (\mmle{k} - \aadapest_{k}) \|^{2r}\\
  &\le & \alpha  C(p,r).
   \end{eqnarray*}
Here \( \EE_{0, \Sigma} \) stands for the expectation w.r.t. \( \norm{0}{\Sigma} \) and \( C(p,r)  \) is given by \eqref{C(p,r)}.
\end{lemma}

\par As before we suppress the dependence on \( x \). To get the propagation property we study for \( k=1, \ldots, K \) the joint distributions of \( \bb e_j^{\T} \mmle{1} , \ldots , \bb e_j^{\T} \mmle{k}\), that is the distribution of \( \bb e_j^{\T} \tilde{\TTa}_k \), the \( j \)th row of the matrix \( \tilde{\TTa}_k \). Obviously,
 \begin{eqnarray*}
   \EE_{\ff, \Sigma_0} [\bb e_j^{\T}\tilde \TTa_k] &=&  \bb e_j^{\T} \TTa^*_k
        =(\bb e_j^{\T}\bbpf{1}, \ldots, \bb e_j^{\T} \bbpf{k}),  \\
   \EE_{\tta, \Sigma} [\bb e_j^{\T} \tilde \TTa_k] &=&  \bb e_j^{\T} \TTa_k
            = (\bb e_j^{\T} \tta, \ldots,\bb e_j^{\T} \tta).
 \end{eqnarray*}
Recall that the matrices \( \SSigma_{k,0} \) and \( \SSigma_{k} \) have a block structure. Now, for instance, to study the estimator of the first coordinate of the vector \( \tta=\tta(x) \), or
of \( f(x) \) in the case of the polynomial basis, we take the first elements of each block and so on.
Denote the \( k \times k \) covariance matrices of the \( j \)th elements of the vectors
\( \mmle{1}, \ldots , \mmle{k} \) by
\begin{eqnarray}
    \SSigma_{k,j}
        &\eqdef& \big\{ \cov_{\tta, \Sigma}\big[\mle{l}{j},\mle{m}{j}   \big] \big\}_{1\le l\le m\le k}\nn
         &=&     \DD_{k,j} (J_k \otimes \Sigma) \DD_{k,j}^{\T}, \label{def SSigma j}\\
    \SSigma_{k,0,j}
        &\eqdef& \big\{ \cov_{\ff, \Sigma_0}\big[\mle{l}{j},\mle{m}{j}   \big] \big\}_{1\le l\le m\le k}\nn
         &=&     \DD_{k,j} (J_k \otimes \Sigma_0) \DD_{k,j}^{\T},\label{def SSigma 0 j}
\end{eqnarray}
where \(J_k \) is a \( k \times k \) matrix with all its elements equal to \(1\), and the \( k \times nk \) block diagonal matrices \( \DD_{k,j} \) is defined by
\begin{eqnarray}
 \DD_{k,j}    &\eqdef&   \bb e_j^{\T} D_1 \oplus \cdots \oplus\bb e_j^{\T} D_k , \label{def DD_K}
            = \big( I_k \otimes \bb e_j^{\T} \big) \DD_k\nn
  D_l     &\eqdef&    \B{l}^{-1} \PPsi \W{l}, \;\;\; l = 1, \ldots ,k.
\end{eqnarray}
Moreover, the following representation holds:
\begin{eqnarray}\label{SSigma_kj via SSigma_k}
  \SSigma_{k,j} &=& \big( I_k \otimes \bb e_j^{\T} \big) \DD_k
                        \big( J_k \otimes \Sigma  \big) \DD_k^{\T}
                            \big( I_k \otimes \bb e_j^{\T} \big)^{\T} \nn
  &=& \big( I_k \otimes \bb e_j \big)^{\T} \SSigma_k \big( I_k \otimes \bb e_j \big),
\end{eqnarray}
where \( \SSigma_k \) is defined by \eqref{def SSigma}. Similarly,
\begin{equation}\label{SSigma_k0j via SSigma_k0}
    \SSigma_{k,0,j} = \big( I_k \otimes \bb e_j \big)^{\T} \SSigma_{k,0} \big( I_k \otimes \bb e_j \big).
\end{equation}
Thus, the important relation \eqref{multi cond sigma} is preserved for \( \SSigma_{k,j} \) and
\( \SSigma_{k,0,j} \) obtained by picking up the \( (j,j) \)th elements of each block of \( \SSigma_{k} \) and \( \SSigma_{k,0} \) respectively.

\par With usual notation \( \gamma^{(j)} \) for the \( j \)th component of \( \gamma \in \RR^k \), denote by
\begin{eqnarray}
  b_j(k) &\eqdef& (\bb e_j^{\T} (\bbpf{1} - \tta), \ldots,\bb e_j^{\T} (\bbpf{k} - \tta) )^{\T} \nn
  &=& ( (\bbpf{1} - \tta)^{(j)}, \ldots, (\bbpf{k} - \tta)^{(j)} )^{\T} \in \RR^k \label{def_bj(k)}  \\
  \Delta_j(k) & \eqdef & b_j(k)^{\T} \SSigma_{k,j}^{-1} \, b_j(k)  \label{def_Delta j(k)}.
\end{eqnarray}
\begin{proposition}\label{Propagation result componentwise theorem} \textbf{``Componentwise'' propagation property}
   \par Under the conditions \( (A1) - (A6) \) and \( (PC) \) for any \( k \le K \) the following upper bound holds:
    \begin{eqnarray}\label{Propagation result component}
      && \left( \frac{nh_{k}^d\Lambda_0}{\bar\s^2_{max}(k) } \right)^{r/2}
               \EE |\bb e_j^{\T}\mmle{k}(x) - \bb e_j^{\T}\aadapest_k(x)|^r\nn
       & \le & (\alpha \EE|\chi^2_p|^r)^{1/2} (1+\delta)^{pk/4}(1-\delta)^{-3pk/4}
        \exp\left\{ \varphi(\delta)
        \frac{\Delta_j(k)}{2(1-\delta)}\right\}
 \end{eqnarray}
 with \( \varphi(\delta) \) as in Theorem \ref{Propagation result theorem}.
 \end{proposition}
\begin{corollary}\label{corr:PP for pol bas}
Let the basis be polynomial. Then under the conditions of the preceding theorem
\(   \EE |\flej{j-1}{k}{x} - \adaplpest^{(j-1)}_k(x) |^r \) satisfy \eqref{Propagation result component}
\end{corollary}

\begin{proof}
The proof essentially follows the line of the proof of Theorem \ref{Propagation result theorem}.
If the distributions of \( \vec \tilde \TTa_k \) were Gaussian, then any subvector is also Gaussian.

\par Denote by
 \begin{equation*}
     \P_{\tta, \Sigma}^{k,j} = \norm{(\bb e_j^{\T}\tta, \ldots, \bb e_j^{\T} \tta)^{\T}}{\SSigma_{k,j}}
\end{equation*}
and by
 \begin{equation*}
     \P_{\ff, \Sigma_0}^{k,j} = \norm{(\bb e_j^{\T}\bbpf{1}, \ldots, \bb e_j^{\T} \bbpf{k} )^{\T}}{\SSigma_{k,0,j}}
\end{equation*}
\( k=1, \ldots, K \), the distributions of \( e_j^{\T} \tilde \TTa_k \) under the parametric assumption and in the non-parametric case.

\par By the Cauchy-Schwarz inequality and Lemma \ref{PC componentwise}
\begin{equation*}
    \left( \frac{nh_{k}^d\Lambda_0}{\bar\s^2_{max}(k) } \right)^{r/2}
               \EE |\bb e_j^{\T}\mmle{k}(x) - \bb e_j^{\T}\aadapest(x)|^r
  \le  (\alpha \EE|\chi^2_p|^r)^{1/2} \big( \EE_{\tta, \Sigma} [Z^2_{k,j}] \big)^{1/2}
\end{equation*}
with the Radon-Nikodym derivative given by \(  Z_{k,j} = \dd \P_{\ff, \Sigma_0}^{k,j}/\dd \P_{\tta, \Sigma}^{k,j} \). By inequalities \eqref{SSigma_kj via SSigma_k} and \eqref{SSigma_k0j via SSigma_k0} the analog of \( (A3) \) is preserved for \( \SSigma_{k,0,j} \) and \( \SSigma_{k,j} \), that is, there exists \( \delta \in [0,1  ) \) such that
\begin{equation}\label{condition Sj}
 (1-\delta)   \SSigma_{k,j}\preceq  \SSigma_{k,0,j} \preceq (1+\delta) \SSigma_{k,j}
\end{equation}
for any \( k \le K \) and \( j=1, \ldots ,p \). Then the assertion of the theorem follows by the Taylor expansion at the point \( (\bb e_j^{\T}\tta, \ldots, \bb e_j^{\T} \tta)^{\T} \) and \eqref{condition Sj} similarly to the proof of Theorem~\ref{Propagation result theorem}.
\end{proof}
At this point we introduce the ``componentwise'' small modeling bias conditions:
\begin{description}
\item[\( \bb{(SMBj)} \)]
\emph{
Let for some \( j = 1, \ldots, p \), \( k(j) \le K \) and  \( \theta^{(j)} = e_j^{\T} \tta \) exist a finite constant \( \Delta_j  \ge 0 \) such that
\begin{equation} \label{smbj}
  \Delta_j(k(j)) \le  \Delta_j,
\end{equation}
where \( \Delta_j(k) \) is defined by \eqref{def_Delta j(k)}.
}
\end{description}
\begin{definition}
For each \( j = 1, \ldots, p \) the oracle index \( k^*(j) \) is defined as the largest index in the scale for which the \( (SMBj) \) condition holds, that is
\begin{equation}\label{def oracle index kj}
    k^*(j) = \max \{ k\le K : \Delta_j(k)\le  \Delta_j \}.
\end{equation}
\end{definition}
\begin{proposition}\label{oracle result componentwise}
Assume \( (A1) - (A6)\) and \( (PC) \).
Let \( h_1 \), the smallest bandwidth, be such that the first estimator \( \bb e_j^{\T}\mmle{1}(x) \) is always accepted by the adaptive procedure.
Let \( k^*(j) \) be the oracle index defined by~\eqref{def oracle index kj}, \( j = 1, \ldots, p \). Then we have for the risk between the \( j \)th coordinates of the adaptive and oracle estimator the following bound:
\begin{eqnarray}\label{oracle ineq compon}
         && \left( \frac{nh_{k^*(j)}^d\Lambda_0}{\bar\s^2_{max}(k^*) } \right)^{r/2}
               \EE |\bb e_j^{\T}\mmle{k^*(j)}(x) - \bb e_j^{\T}\aadapest(x)|^r\\
       & \le &   \z_{k^*(j)}^{r/2} +
       (\alpha \EE|\chi^2_p|^r)^{1/2} (1+\delta)^{pk_j^*/4}(1-\delta)^{-3pk_j^*/4}
        \exp\left\{ \varphi(\delta)
        \frac{\Delta_j}{2(1-\delta)}\right\}, \nonumber
\end{eqnarray}
where \( \varphi(\delta) \) is as in Theorem \ref{Propagation result componentwise theorem}.
\end{proposition}
\begin{corollary}\label{corr:oracle risk bound for pol bas}
Let the basis be polynomial. Then under the conditions of the preceding theorem the risk between the adaptive and oracle estimators
\begin{equation*}
     \EE |\flej{j-1}{k^*(j)}{x} - \adaplpest^{(j-1)}(x)|^r
\end{equation*}
 satisfy \eqref{oracle ineq compon}.
\end{corollary}
\begin{remark}
 The statements of this and the preceding proposition are of the same type that their vector counterparts. They are needed for asymptotical results of the last section.
 \end{remark}
\begin{proof} To simplify the notation we suppress the dependence on \( j \) in the index \( k \). Similarly to the proof of Theorem \ref{oracle result} we consider the disjunct events
\( \{ \adapind \le k^*\} \) and \( \{ \adapind > k^*\} \). Therefore,
\begin{eqnarray*}
       && \EE|\bb e_j^{\T}\mmle{k^*}(x) - \bb e_j^{\T}\aadapest(x)|^r\\
&=& \EE|\bb e_j^{\T}\mmle{k^*}(x) - \bb e_j^{\T}\aadapest(x)|^r \,\ind\{ \adapind \le k^* \} \\
&+&
\EE|\bb e_j^{\T}\mmle{k^*}(x) - \bb e_j^{\T}\aadapest(x)|^r \, \ind \{ \adapind > k^*\}.
\end{eqnarray*}
By Lemma \ref{bound for the componentwise differences} and the definition of the test statistic \( T_{k^*, \adapind} \) the second summand can be easily bounded:
\begin{eqnarray*}
    && \left( \frac{nh_{k^*}^d\Lambda_0}{\bar\s^2_{max}(k^*) } \right)^{r/2}
          \EE|\bb e_j^{\T}\mmle{k^*}(x) - \bb e_j^{\T}\aadapest(x)|^r \,\ind \{ \adapind > k^*\} \\
    &\le & \EE \| \B{k^*}^{1/2} (\mmle{k^*}(x) - \aadapest(x)) \|^r \,\ind \{ \adapind > k^*\}\\
    &\le &  \z_{k^*}^{r/2}.
\end{eqnarray*}
 To bound the first summand we use the ``componentwise'' analog of Theorem~\ref{Propagation result theorem}, namely Theorem~\ref{Propagation result componentwise theorem} that completes the proof.
\end{proof}
\subsection{SMB and the bias-variance trade-off}
\label{subsect:SMB and the bias-variance trade-off}
It was shown in~\cite{SV} that the small modeling bias (\( SMB1 \) here) condition given by~\eqref{smbj} can be obtained from the ``bias-variance trade-off'' relations. Notice that our set-up includes the set-up from \cite{SV} as a particular
 (\( p=0 \), \( \delta =0 \), \( \s(\cdot) \equiv \s \) is a known constant) case. To prove that the similar relation holds in the present case we need the following definition. Let the basis be polynomial. Given a point \( x \) and method of localization \( w \), for any \( j = 1, \ldots, p \) the ``ideal adaptive bandwidths'', see \cite{LepMamSpok97} and \cite{LepSpok97}, is defined as follows:
\begin{equation}\label{kj star}
    k^\star(j) = \max\{ k\le K : \bar{b}_{k,f^{(j-1)}}(x) \le  C_j(w) \s_{k,j}(x) \sqrt{d(n)} \},
\end{equation}
where \( C_j(w) \) is a constant depending on the choice of the smoother \( w \),
\begin{eqnarray*}
  \bar{b}_{k,f^{(j-1)}}(x) &=& \sup_{1\le l\le k}|\bb e_j^{\T} \bbpf{l}(x) -  f^{(j-1)}(x)| , \\
  \s_{k,j}^2(x) &=& \Var_{\ff, \Sigma_0} [\bb e_j^{\T} \mmle{k}(x) ],\\
  d(n) &=& \log(h_K/h_1),
\end{eqnarray*}
and \( f^{(0)} \) stands for the function \( f \) itself. To bound the ``modeling bias'' \( \Delta_j(k) \) we need the following assumption:
\begin{description}
\item[\( \bb{(A7)} \)]
\emph{
There exists  a constant \( s_j > 0 \) such that for all \( k \le K \)
\begin{equation}\label{sigma kj}
     \SSigma^{-1}_{k,j} \preceq s_j \SSigma^{-1}_{k,j,diag}
\end{equation}
}
\end{description}
where \( \SSigma_{k,j,diag} = \diag\big(\Var_{\tta, \Sigma}[\bb e_j^{\T} \mmle{1}(x)] ,
 \ldots , \Var_{\tta, \Sigma}[\bb e_j^{\T} \mmle{k}(x)] \big) \) is a diagonal matrix composed of the diagonal elements of \( \SSigma_{k,j} \).
\begin{remark}
In order to understand the meaning and fulfillment of this assumption let us consider for simplicity the case of
\( p=0 \) (local constant fitting). Then \eqref{sigma kj} can be rewritten as
\begin{equation*}
    \exists \, s>0: \, \bb R_k = \SSigma_{k,diag}^{-1/2} \SSigma_{k} \SSigma_{k,diag}^{-1/2}
                                                        \succeq s^{-1} \bb I_k \;\;\; \forall k \le K,
\end{equation*}
where \( \SSigma_{k} = (\cov [\mmle{l}, \mmle{m}])_{1 \le l \le m \le k} \) is a \( k \times k \) positive definite matrix, \( \SSigma_{k,diag} = \diag(v_1, \ldots , v_k) \) with \( v_l = \Var_{\tta, \SSigma}[\mmle{l}] > 0\), \( l= 1, \ldots, k \). We immediately see the following:
\begin{enumerate}
  \item Since \( \bb R_k \succ 0 \), it is known that for any symmetric matrix \( \bb A \) one can find a sufficiently small in absolute value real number \( \tau \) s.t. \( \bb R_k - \tau \bb A \succ 0 \).
  \item \( \bb R_k = (\rho_{lm})_{1 \le l \le m \le k} \) is a correlation matrix with entries
  \begin{equation*}
    \rho_{lm} = (v_l v_m)^{-1/2} \cov [\mmle{l}, \mmle{m}].
  \end{equation*}
  Moreover, \( 1 \ge \rho_{lm} > 0  \) since for \( \w{k}{i} \in [0,1] \) the estimators are strictly positively correlated. Indeed,
       \begin{eqnarray*}
              \mmle{k} & = &(\sum_i \frac{\w{k}{i} } {\s_i^2})^{-1} \sum_i \frac{\w{k}{i} } {\s_i^2} \Yi \\
         \cov [\mmle{l}, \mmle{m}] &=& (\sum_i \frac{\w{l}{i} } {\s_i^2})^{-1}(\sum_i \frac{\w{m}{i} } {\s_i^2})^{-1} \sum_i \frac{\w{l}{i}\w{m}{i} \s_{0,i}^2} {\s_i^4} > 0.
       \end{eqnarray*}
Strict inequality takes place because the estimators have a common support and therefore are dependent. Below we shall see that \( (A7) \) essentially means that the estimators should not be correlated too strongly, which in its turn is provided by the assumption on the ``geometrical growth of the scales'', i.e. by \( (A4) \). Indeed, since \( \rho_{lm} > 0 \), we have by direct calculations
\begin{equation*}
    (1-\rho_{max}) \bb I_k  \preceq \bb R_k \preceq (1-\rho_{max}) \bb I_k  + \rho_{max} \bb J_k,
\end{equation*}
where \( \rho_{max} = \max_{1 \le l < m \le k}\{ \rho_{lm} \} \) is the maximal correlation of the off-diagonal elements of \( \SSigma_{k} \) and \( \bb J_k \) is a \( k \times k \) matrix with all its elements equal to one. Thus we see that \( s = (1-\rho_{max})^{-1} \) explodes when the maximal correlation (except for the variations) is close to one.
  \item Connection with \( (A4) \). Assume additionally that the weights \( \w{l}{i} = \ind \{ \| \Xi -x \| \le h_l/2 \} \). Then \( \rho_{lm}^2 = \sum_i \w{l}{i} /(\sum_i \w{m}{i})  \) for \( l <m \) and \( v_m/v_l = \B{l}/\B{m} = \rho_{lm}^2 \).
       Since \( u \ge u_0>1 \), assumption \( (A4) \) provides \( 0 < u^{-(m-l)/2 }\le \rho_{lm} \le  u_0^{-(m-l)/2 } < 1 \) .
\end{enumerate}
\end{remark}
We have the following result:
\begin{lemma}\label{th:SMB from the balance equation}
  Let the weights \( \{ \w{k}{i}(x) \} \) satisfy~\eqref{binar product} and the basis be polynomial \( \{ 1,t-x,(t-x)^2/2!, \ldots, (t-x)^{p-1}/(p-1)! \} \).
  Granted assumptions \( (A1) \) -- \( (A4) \) and \( (A7) \) for any (possibly fixed) \( n \), any given point \( x \), smoothing function \( w \) and \( j = 1, \ldots, p \) the choice of \( k(j) = k^\star(j) \) defined by \eqref{kj star} with \( d(n)=1 \) implies the \( (SMBj) \)
condition \( \Delta_j(k(j)) \le \Delta_j \) with
the constant \( \Delta_j < 2 s_j C^2_j(w) (1-u_0^{-1})^{-1}<\infty \).
\end{lemma}
\begin{proof}
Consider the quantity \( b_j(k)^{\T} \SSigma^{-1}_{k,j,diag} b_j(k)\). For the polynomial basis \( \bb e_j^{\T} \tta(x) = f^{(j-1)} (x) \). In view of \eqref{binar product} the matrix
\( \SSigma_{k,j,diag} \) is particularly simple:
\begin{eqnarray*}
  \SSigma_{k,j,diag}
    &=&
        \diag (\bb e_j^{\T} \B{1}^{-1} \bb e_j , \ldots,  \bb e_j^{\T} \B{k}^{-1} \bb e_j) \\
    &=& \diag(\Var_{\tta, \Sigma}[\mmle{1}^{(j)}(x)] , \ldots, \Var_{\tta, \Sigma}[\mmle{k}^{(j)}(x)] ),
\end{eqnarray*}
that is \( \SSigma_{k,j,diag} \) is a diagonal matrix of the variances of the \( j \)th coordinates of vectors \( \mmle{1} , \ldots, \mmle{k}\). Then by~\( (A4) \) and \eqref{Vk bound}
\begin{eqnarray*}
  b_j(k)^{\T} \SSigma^{-1}_{k,j,diag} b_j(k)
        &=&
            \sum_{l=1}^k   \frac{|\bb e_j^{\T}\bbpf{l} - f^{(j-1)}(x)|^2}{\bb e_j^{\T} \B{l}^{-1} \bb e_j}
        \\
        &\le &
            \big(\bar{b}_{k,f^{(j-1)}}(x)\big)^2
                \sum_{l=1}^k   \frac{1}{\bb e_j^{\T} \B{l}^{-1} \bb e_j}   \\
        &\le & \frac{\big(\bar{b}_{k,f^{(j-1)}}(x)\big)^2}{\bb e_j^{\T} \B{k}^{-1} \bb e_j}
                \sum_{l=1}^k u_0^{-(k-l)}\\
        &\le & \frac{\big(\bar{b}_{k,f^{(j-1)}}(x)\big)^2 (1+\delta)}{\s_{k,j}^2(x) (1-u_0^{-1})}.
\end{eqnarray*}
By \eqref{kj star} with \( d(n) =1 \) the choice of \( k=k^\star(j) \)
implies \( \big(\bar{b}_{k,f^{(j-1)}}(x)\big)^2 \le  C^2_j(w) \s_{k,j}^2(x)\). Thus
\begin{equation*}
     b_j(k)^{\T} \SSigma^{-1}_{k,j,diag} b_j(k) \le  (1+\delta) C^2_j(w)(1-u_0^{-1})^{-1}
\end{equation*}
and
\begin{eqnarray*}
  \Delta_j(k) &=&  b_j(k)^{\T} \SSigma^{-1}_{k,j} b_j(k) \le s_j C^2_j(w) (1+\delta) (1-u_0^{-1})^{-1} \\
   &<&  2 s_j C^2_j(w) (1-u_0^{-1})^{-1} < \infty,
\end{eqnarray*}
since \( u_0 \) from \( (A4) \) is strictly larger than \( 1 \).
\end{proof}
\begin{remark}\label{Rem_SMB and indicator kernels}
The assumption that the weights \( \{ \w{k}{i}(x) \} \) satisfy~\eqref{binar product}, that is that they are of the indicator-type, seems to bee too restrictive. This assumption allows to show the connection between the small modeling bias condition and the classical bias-variance trade-off without technical complications for any~\( n \), including the case of the fixed sample size. Relaxing of the consideration to the asymptotic case does not require such an assumption on the weights, see the lemmas below. Moreover, since this section essentially serves for checking the rate of convergence of the adaptive estimator at a point w.r.t. the H\"older classes of functions and since by \( (A2) \) the windows are nested, to get the first impression it is enough to consider the design in \( \RR \), as in the case of the nested windows the generalization of the adaptive procedure to~\( \RRd \) is straightforward. On the contrary non-nested windows that are related to estimation on anisotropic classes require drastic modifications of the procedure, see~\cite{KLP2001} and~\cite{KLP2007}.
\end{remark}
\begin{lemma}
Let the basis be polynomial and for each \( k \) the weight function \( \w{k}{\cdot} (x) = W((\cdot - x) h_k^{-1})\)
be non-negative, bounded with \( \supp W(\cdot) \subset [0,1] \) and such that the Lebesgue measure of the set \( \{ u: W(u)^2 >0 \} \) is strictly positive. Let \( \Xi = i/n \), \( i=1, \ldots, n \), and \( h_k = h_k(n) \) be a sequence s.t. \( h_k(n) \to 0 \) and \(n h_k(n) \to \infty\) as \( n \to \infty \). Let the variance be either known \( (\s_i \equiv \s_{0,i} = \s(\cdot)) \) and continuous at the neighborhood of \( x \), either the known ``model'' variance be locally bounded: i.e. \( \exists \, 0< \s_{min}(k) \le  \s_{max}(k) < \infty \) s.t. \( \s_{min}(k) \le \s_i \le \s_{max}(k) \) for \( \forall \, i:\w{k}{i} (x) > 0\). For a square matrix \( A \) by \( A_{diag} \) we denote a diagonal matrix with the same entries as the main diagonal of \( A \). Then
\begin{enumerate}
  \item
  \begin{equation*}
    \bb e_j^{\T} \Var[\mmle{k} (x)] \bb e_l =  O\Big( \frac{\s^2(x)}{n h_k^{j+l-1}} \Big)
                                            =O\Big( \bb e_j^{\T} \B{k}^{-1} \bb e_l \Big),
  \end{equation*}
  as \( n \to \infty \);
  \item For \( n \) sufficiently large we have
  \begin{eqnarray*}
      && \s_{max}^{-2}(k) \diag(\mu_1(W), h_k^{2}\mu_2(W), \ldots , h_k^{2(p-1)}\mu_{2(p-1)}(W)) \\
    &\lesssim& (nh_k)^{-1}(\B{k} (x))_{diag}\\
     &\lesssim & \s_{min}^{-2}(k) \diag(\mu_1(W), h_k^{2}\mu_2(W), \ldots , h_k^{2(p-1)}\mu_{2(p-1)}(W))
  \end{eqnarray*}
  with the moments of the kernel \( W(\cdot) \) defined by
  \begin{equation*}
    \mu_{\pi}(W) = \int u^{\pi}W(u) \dd u;
  \end{equation*}
    \item By \eqref{Vk def} \( \Var \mmle{k} = \B{k}^{-1} \tilde{\B{k}}  \B{k}^{-1} \), where
    \( \tilde{\B{k}} =\PPsi\W{k} \Sigma_0 \W{k}  \PPsi^{\T} \) is a Gram matrix
    (c.f. Remark~\ref{Rem_Gram matrix}) and therefore the H\"older inequality is applicable to its off-diagonal elements. Since
    \begin{equation*}
        (\tilde{\B{k}})_{diag}
        =
        \diag \Big(\sum_{i=1}^n \frac{\w{k}{i}^2 (x) }{\s_i^2} \frac{{\s_{0,i}^2}}{{\s_i^2}}, \ldots ,
        \sum_{i=1}^n \frac{(\Xi -x)^{2(p-1)}}{((p-1)!)^2} \frac{\w{k}{i}^2 (x) }{\s_i^2} \frac{{\s_{0,i}^2}}{{\s_i^2}}
        \Big)
    \end{equation*}
    and assuming \( (A3) \) similarly to the statement 2 we have for \( n \) sufficiently large
    \begin{equation*}
        (nh_k)^{-1} \bb e_j^{\T}\tilde{\B{k}}\bb e_j \lesssim \frac{1+\delta}{\s_{min}^{2}(k)}  \diag(\mu_1(W^2), h_k^{2}\mu_2(W^2), \ldots , h_k^{2(p-1)}\mu_{2(p-1)}(W^2))
    \end{equation*}
    and the bounds for the variance of \( j \)th coordinate of \( \mmle{k} \):
    \begin{equation*}
  \frac{(1-\delta) \s_{min}^{2}(k)}{nh_k^{1+2(j-1)}}   \lesssim \bb e_j^{\T} \Var [\mmle{k}]\bb e_j
            \lesssim \frac{(1+ \delta) \s_{max}^{2}(k)}{nh_k^{1+2(j-1)}}   .
    \end{equation*}
That is \( \bb e_j^{\T} \Var [\mmle{k}]\bb e_j = O\Big( \bb e_j^{\T} \B{k}^{-1}  \bb e_j\Big) \). The constants depend on \( \s_{min}^{2}(k) \), \( \s_{max}^{2}(k) \) and the moments of \( W \) and \( W^2 \).
\end{enumerate}
\end{lemma}
\begin{remark}\label{Rem_substit of var by Bk}
When the constants are not the target in the study of rate of convergence the last display allows to substitute in the balance equation \eqref{kj star} the variance by the \( (j,j) \)th component of \( \B{k}^{-1} \), with the proviso that \( \delta \) is ``well behaved'', c.f. Remark~\ref{Rem_relative_noise_error}.
\end{remark}
\begin{proof}
The statement of the lemma and its proof is essentially in the spirit of the Theorem 2.1 in \cite{RuppertWand} and Theorem 3.1 in \cite{Fan and Gijbels book}, where the study was performed for the random design.
\end{proof}
\begin{lemma}\label{Lm_converg to B}
Let for each \( k \) the weight function \( \w{k}{\cdot} (x) = W((\cdot - x) h_k^{-1})\)
be non-negative, bounded with \( \supp W (\cdot) \subset [0,1] \) and such that the Lebesgue measure of the set
\( \{ u: W(u) >0 \} \) is strictly positive. Let \( \Xi = i/n \), \( i=1, \ldots, n \), and \( h_k = h_k(n) \) be a sequence s.t. \( h_k(n) \to 0 \) and \(n h_k(n) \to \infty\) as \( n \to \infty \). Let \( \Psi(u) = (1,u,..., u^{p-1}/(p-1)! )^{\T} \) and \( \Psii \eqdef \Psi(i/n - x) \).
\begin{enumerate}
  \item Denote by
\( \B{k}^{\sharp} = \B{k}^{\sharp} (x) = \PPsi  \cc W_{k}\PPsi^{\T} =
        \sum_{i=1}^{n} \Psii \Psii^{\T} \w{k}{i}\). Then with \( H = \diag(1,h_k,\ldots,h_k^{p-1}) \) we have
  \begin{equation*}
    (nh_k)^{-1} H^{-1} \B{k}^{\sharp} H^{-1} \to \bb B=\int \Psi(u)\Psi^{\T}(u)W(u) \dd u
  \end{equation*}
 as \( n \to \infty \) , where the matrix \( \bb B \) is positive definite and independent on \( x \) and \( n \).
  \item Moreover, assuming the known ``model'' variance be locally bounded: i.e. \( \exists \, 0< \s_{min}(k) \le  \s_{max}(k) < \infty \) s.t.\( \s_{min}(k) \le \s_i \le \s_{max}(k) \) for \( \forall i:\w{k}{i} (x) > 0\) we have for sufficiently large \( n \):
      \begin{equation*}
        0\prec \s_{max}^{-2}(k) \bb B \preceq  (nh_k)^{-1} H^{-1} \B{k} H^{-1} \preceq \s_{min}^{-2}(k) \bb B.
      \end{equation*}
\end{enumerate}
\end{lemma}
\begin{proof}
The first statement of the lemma is based on the convergence of Riemann sums. The non-degenerateness of \( \bb B \) is the Lemma 1.4 in \cite{Tsybakov} and follows from the fact that the polynomials of degree \( \le p-1 \) have at most \( p-1 \) different zeros.
\par To justify the second statement it is enough to remark that \( \s_{max}^{-2}(k) \B{k}^{\sharp} \preceq \B{k} \preceq \s_{min}^{-2}(k) \B{k}^{\sharp}\) and that the first statement implies \( (nh_k)^{-1} \gamma^{\T} H^{-1} \B{k}^{\sharp} H^{-1} \gamma \to \gamma^{\T} \bb B \gamma \) for any vector \( \gamma \).
\end{proof}

\begin{remark}
Using the standard technique it is easy to derive from the above result that for estimation of functions over H\"older classes the methodology proposed in \cite{KatkSpok} and \cite{SV} and generalized in the present paper delivers the minimax rate of convergence up to a logarithmic factor, see the following subsection for details.
\end{remark}
\subsection{Rates of convergence}
At this section \( d=1 \) and the basis is polynomial with the columns of the design matrix \( \PPsi \) given by
\begin{equation*}
     \Psii= \Psi(\Xi-x) = \left(1,\, \Xi -x, \ldots, (\Xi -x)^{\p-1}/(p-1)!\right)^{\T}.
\end{equation*}.
The polynomial weights \( W^*_{l,\,i} \) are given by
\begin{equation}\label{polynomial weights_multiscaled}
    W^*_{l,\,i}(x) = \bb e^{\T}_1 \B{l}^{-1} \Psii  \w{l}{i}(x)/\s^2_i
\end{equation}
with \( \B{l}\) defined by \eqref{B} and the variance term given by
\(     \s^2_l(x) \eqdef \EE_f[|\bb e^{\T}_1 \mmle{l}(x)  - \bb e^{\T}_1 \bbpf{l}(x) |^2] \). Here
\begin{equation*}
    \bb e^{\T}_1 \bbpf{l}(x) = \EE_f [\fle{l}{x}] = \sum_{i=1}^{n}  W^*_{l,\,i}(x) \ffi
\end{equation*}
is a local linear smoother of the function \( f  \) at the point \( x \) corresponding to \( l \)th scale. Define the ``monotonized'' bias by
\begin{equation}\label{adaptive bias}
    \bar{b}_{k,f}(x) 
                     =      \sup_{1\le l \le k } |\bb e^{\T}_1 \bbpf{l}(x) - f(x)|
\end{equation}

\par Before proceeding with analysis of the convergence rate we need to derive bounds for the bias and variance.

\begin{description}
\item[\( \bb{(A8)}\)]
Let the known ``model'' variance be locally bounded: i.e. \( \exists \, 0< \s_{min}(k) \le  \s_{max}(k) < \infty \) s.t. \( \s_{min}(k) \le \s_i \le \s_{max}(k) \) for \( \forall \,  i:\w{k}{i} (x) > 0\).
\item[\( \bb{(A9)} \)]
\emph{There exists a real number \( a_0>0 \) such that for any interval \( A \subseteq [0,1] \) and all \( n \ge 1 \)}
\begin{equation*}
    \frac{1}{n} \sum_{i=1}^n \ind \{ \Xi \in A \} \le a_0 \max \big\{ \int_A\dd t, \frac{1}{n} \big\}.
\end{equation*}
\item[\( \bb{(A10)} \)]
\emph{ The localizing functions (kernels) \( \w{k}{i} \) have compact supports belonging to \( [0,1] \):}
\begin{equation*}
    \w{k}{i}(x) = 0 \;\;\;\text{if} \;\;\; |\Xi - x|> h_k.
\end{equation*}
This immediately implies the similar property for the local polynomial weights:
\begin{equation*}
    W^*_{k,i}(x) = 0 \;\;\;\text{if} \;\;\; |\Xi - x|> h_k.
\end{equation*}
\item[\( (A11) \)]
\emph{ There exists a finite number \( w_{max} \) such that}
\begin{equation*}
    \sup_{k,i}|\w{k}{i}(x) | \le w_{max}.
\end{equation*}

\end{description}
\begin{remark}
Assumption \( (A3) \) implies that the conditional number
\begin{equation}\label{kond_number}
    \kappa(\Sigma) \eqdef \frac{\s^2_{max}}{\s^2_{min}}
\end{equation}
of covariance matrix from the misspecified model \eqref{PA} is finite.
\end{remark}
\begin{lemma}\label{nonuniform upper bounds for bias and variance}
 Assume \( (A1) - (A3) \), \( (A6) \) and \( (A8)-(A11) \). Let \( h'_1 \) be the smallest bandwidth providing \( (A6) \) and \( h''_1 \) be the smallest bandwidth s.t. the first estimator \( \mmle{1} \) is accepted by the adaptive procedure. Denote by  \( h_1 \ge \max\{  1/(2n), h'_1, h''_1 \} \).
 Let the regression function \( f(\cdot) \) belong to the H\"older class
  \(  \Sigma (\beta, L) \) on \( [0,1] \), and let \( \{ \fle{k}{x} \}_{k=1}^K\) be the \( LP_k(p-1) \) estimators of \( f(x) \) with \( p -1 = \lfloor \beta \rfloor\). Then for sufficiently large \( n \), and any \( h_k \) satisfying \( h_K > \ldots > h_k > \ldots > h_1 \), \( k=1, \ldots, K \), we have
\begin{eqnarray*}
  |\bar{b}_{k,f}(x)|    &\le & C_2  \kappa(\Sigma) \frac{Lh_k^{\beta}}{(p-1)!}, \\
  \s^2_k(x)   &\le & (1+\delta) \frac{\s^2_{max} }{nh_k \Lambda_0}
\end{eqnarray*}
with  \( C_2 = 2 w_{max} a_0 \sqrt{e}/\Lambda_0 \) and \( \delta \in [0,1) \) from \( (A3) \).
\end{lemma}
The proof is moved to Appendix.
\begin{proposition}\label{Proposition_rate in lin mod}
Let the model \eqref{PA} be satisfied.  Assume \( (A1) - (A4) \), \( (A6) \) and \( (A8)-(A11) \). Let \( h'_1 \) be the smallest bandwidth providing \( (A6) \) and \( h''_1 \) be the smallest bandwidth s.t. the first estimator \( \mmle{1} \) is accepted by the adaptive procedure. Denote by  \( h_1 \ge \max\{  1/(2n), h'_1, h''_1 \} \).
 Let the regression function \( f(\cdot) \) belong to the H\"older class
  \(  \Sigma (\beta, L) \) on \( [0,1] \), and let \( \{ \fle{k}{x} \}_{k=1}^K\) be the \( LP_k(p-1) \) estimators of \( f(x) \) with \( p -1 = \lfloor \beta \rfloor\). Then for sufficiently large \( n \) for adaptive estimator obtained by the procedure we have
  \begin{equation*}
    \EE|\hat f(x)  - f(x)|^r \asymp \Big(\frac{\log n}{n } \Big)^{\frac{\beta r}{2 \beta +1} }.
  \end{equation*}

\end{proposition}
\begin{proof}
If the model \eqref{PA} is true, then \( \Delta = 0 \) and one can take \( k^* \) from \eqref{kj star} with \( d(n) = \log n \) leading in view of the preceding lemma to the choice of the optimal bandwidth \( h_{k^*}(x) \)
of order \( (\log n /n)^{\frac{1}{2 \beta +1}} \). The oracle bound of Proposition~\ref{oracle result componentwise} gives
\begin{equation*}
     \EE |\hat f(x)  - \tilde f_{k^*}(x)|^r \lesssim \Big(\frac{\log n}{n h_{k^* (x)}} \Big)^{r/2}.
\end{equation*}

Since \( k^* \) is an unknown but deterministic \( \tilde f_{k^*}(x) \) is a standard local polynomial estimator, therefore its quality of estimation is known:
\begin{equation*}
    \EE |f(x)  - \tilde f_{k^*}(x)|^r
        \lesssim (1 /n)^{\frac{\beta r}{2 \beta +1}}\ll (\log n /n)^{\frac{\beta r}{2 \beta +1}}
\end{equation*}
and the assertion follows by application of \( (a+b)^r \le C(r)(a^r + b^r)  \), \( a,b \ge 0 \), \( C(r) = 2^{r-1} \) for \( r\ge 1 \), and equal to one for \( 0<r<1 \), since the rate \( \Big(\frac{\log n}{n } \Big)^{\frac{\beta r}{2 \beta +1} } \) is known to be optimal, c.f. \cite{Lep1990}.
\end{proof}

\begin{proposition}\label{Proposition_rate in nonparametric mod}
 Assume \( (A1) - (A4) \), \( (A6) -(A11) \) and \( \delta = O(1/\log n) \). Let \( h'_1 \) be the smallest bandwidth providing \( (A6) \) and \( h''_1 \) be the smallest bandwidth s.t. the first estimator \( \mmle{1} \) is accepted by the adaptive procedure. Denote by  \( h_1 \ge \max\{  1/(2n), h'_1, h''_1 \} \).
 Let the regression function \( f(\cdot) \) belong to the H\"older class
  \(  \Sigma (\beta, L) \) on \( [0,1] \), and let \( \{ \fle{k}{x} \}_{k=1}^K\) be the \( LP_k(p-1) \) estimators of \( f(x) \) with \( p -1 = \lfloor \beta \rfloor\). Then for sufficiently large \( n \) for the adaptive estimator delivered by the procedure we have
  \begin{equation*}
    \EE|\hat f(x)  - f(x)|^r \lesssim  \Big(\frac{\log^{\gamma} n}{n } \Big)^{\frac{\beta r}{2 \beta +1} }
  \end{equation*}
with \( \gamma = (2 \beta + 1)/(2\beta) \).
\end{proposition}
\begin{proof}
Since now we need to have (SMB) fulfilled, we have to take \( k^* \) from \eqref{kj star} with \( d(n) = 1 \) leading to the suboptimal choice of  \( h_{k^*}(x) \) of order \( (1/n)^{\frac{1}{2 \beta +1}} \) and the assertion follows.
\end{proof}


\

\section{Appendix}\label{section:Auxiliary results}

\subsection{Pivotality and local parametric risk bounds}
\begin{lemma}\label{Pivotality property} {\it \textbf{Pivotality property}}
\par Let \( {(A2)} \) hold. Let \( \bbpf{1} = \cdots = \bbpf{\kappa} = \tta \) for \( \kappa \le K \). Then for any \( k \le   \kappa \) the risk associated with the adaptive estimate at every step of the procedure does not depend on the parameter \( \tta \):
\begin{eqnarray*}
     \EE_{\tta} | ( \mmle{k} -\aadapest_{k} )^{\T}
    \B{k}
        ( \mmle{k} -\aadapest_{k} ) |^{r}
    &=& \EE_{0} | ( \mmle{k} -\aadapest_{k} )^{\T}
    \B{k}
    ( \mmle{k} -\aadapest_{k} ) |^{r},
\end{eqnarray*}
where \( \EE_0 \) denotes the expectation w.r.t. the centered measures \( \norm{0}{\Sigma} \) or \( \norm{0}{\Sigma_0} \).
\end{lemma}
\begin{proof}
At each step \( k \) of the procedure the adaptive estimator \( \aadapest_{k} \) coincides with one of the nonadaptive estimators \( \mmle{1}, \ldots, \mmle{k} \). If \( \aadapest_{k} = \mmle{k} \), this means that the deviation from the parametric model is not significant and the procedure passes to the next step. On the contrary,
\( \aadapest_{k} = \mmle{m} \) for \( m < k \) means that for some \( l \le m \) the value of the test statistic
\( T_{l,\,m+1}   \) is strictly larger than the threshold \( \z_l \) and the procedure had terminated. Thus one can write the following decomposition:
\begin{eqnarray*}
&& \EE_{\tta} |(\mmle{k} - \aadapest_k )^{\T} \B{k} (\mmle{k} - \aadapest_k ) |^r \\
    &=&
    \sum_{m=1}^{k} \EE_{\tta} \|  \B{k}^{1/2} (\mmle{k} - \mmle{m} ) \|^{2r} \I{\{ \aadapest_{k} = \mmle{m}\}}\\
    &=&\sum_{m=1}^{k-1} \EE_{\tta} \|  \B{k}^{1/2} (\mmle{k} - \mmle{m} ) \|^{2r}
        \I{\{ \exists \, l \le m : \|  \B{k}^{1/2} (\mmle{k} - \mmle{m} ) \|^{2} > \z_l\} }.
\end{eqnarray*}
In the last line the definition of \( T_{l,\,m+1}   \) given by \eqref{Tlk} is used. Since for any \( k \le \kappa \) under the assumptions of lemma
\( \mmle{k} = \tta + \B{k}^{-1} \PPsi \W{k} \Sigma_0^{1/2} \eeps \), the value of \( \tta \) cancels in the differences \( \mmle{k} - \mmle{m} \) and \( \mmle{l} - \mmle{m+1} \) for all \( l \le m < k \), and therefore can be taken equal to zero.
\end{proof}

\par To justify the statistical properties of the considered procedure we need the following simple observation. Let for any \( \tta \), \( \tta' \in \Theta\) the corresponding log-likelihood ratio \( \LL(\W{k},\tta, \tta') \) be defined by \eqref{log-likelihood ratio}.
       Then
\begin{equation*}
    2 \LL(\W{k},\tta, \tta')=
                \| \W{k}^{1/2} ( \YY -\PPsi^{\T} \tta'  ) \|^2
     -   \| \W{k}^{1/2}( \YY -\PPsi^{\T} \tta  ) \|^2.
     \end{equation*}
\begin{lemma}\label{Th. Spokoiny_fitted likelihood}
\textbf{Quadratic shape of the fitted log-likelihood}

\par Let for every \( k = 1, \ldots, K\) the fitted log likelihood (FLL) be defined as follows:
\begin{equation*}
    \LL(\W{k},\mmle{k}, \tta') \eqdef \max_{\tta \in \Theta} \LL(\W{k},\tta, \tta').
  \end{equation*}
Then
\begin{equation}
    2 \LL(\W{k},\mmle{k}, \tta)
    =  ( \mmle{k} -\tta )^{\T} \B{k} ( \mmle{k} -\tta ).
\end{equation}
\end{lemma}

\begin{proof}
Notice that \( \LL(\W{k}, \tta) \) defined by \eqref{log-likelihood} is quadratic in \( \tta \). The assertion follows from the second order Taylor series expansion around the point \( \mmle{k} \), because it is the point of maximum, and the second derivative is the constant matrix \( \B{k} \).
\end{proof}
Let the matrix \( \S \) be defined as follows:
\begin{equation}\label{def S}
    \S \eqdef \Sigma_0^{1/2} \W{k} \PPsi^{\T} \B{k}^{-1}
        \PPsi \W{k} \Sigma_0^{1/2} .
\end{equation}
Then for the distribution of \( \LL(\W{k},\mmle{k}, \bbpf{k}) \) one observes so-called ``Wilks phenomenon'', c.f.~\cite{FanZhangZhang}, described by the following theorem:

\begin{proposition}\label{Wilks theorem}
Let the regression model be given by \eqref{true model} and the parameter \( \bbpf{k} = \bbpf{k}(x) \) maximizing the expected local log-likelihood be defined by \eqref{prameter of BPF}. Then for any \( k = 1, \ldots, K \) the following equality in distribution takes place:
\begin{equation}\label{likelihood = weightsum of sqnormal }
    2 \LL ( \W{k}, \mmle{k}, \bbpf{k} ) \eqdistr
    \lambda_{1}(\S)\bar\varepsilon_{1}^{2} + \cdots +
    \lambda_{p}(\S)\bar\varepsilon_{p}^{2}
\end{equation}
with \( p = \rank (\B{k}) = \dim \Theta = \p \). Here \( \lambda_{1}(\S) , \ldots , \lambda_{p}(\S) \)
are the non-zero eigenvalues of the matrix \( \S  \), and \( \bar\varepsilon_{i} \) are independent
standard normal random variables.

\par Moreover, under \( (A3) \) the maximal eigenvalue \( \lambda_{max}(\S) \le 1+ \delta \), and for any \( \z > 0 \)
\begin{equation}\label{chi squared domination}
    \P \left\{2 \LL ( \W{k}, \mmle{k}, \bbpf{k} ) \ge  \z\right\} \le
    \P \left\{\eta \ge  \z/(1+ \delta)\right\},
\end{equation}
where \( \eta  \) is a random variable distributed according to the \( \chi^2 \)
law with \( p \) degrees of freedom.
\end{proposition}
\begin{remark}
Generally, if \( \B{k} \) is degenerated, the number of terms in \eqref{likelihood = weightsum of sqnormal } is \( p \le \dim \Theta \).
\end{remark}

\begin{proof}
    By Lemma~\ref{Th. Spokoiny_fitted likelihood} and the decomposition \eqref{linearity if quasiMLE}
     it holds that:
\begin{eqnarray*}
 2 \LL ( \W{k}, \mmle{k}, \bbpf{k} )
 &=& ( \mmle{k} -\bbpf{k} )^{\T} \B{k}
        ( \mmle{k} -\bbpf{k} )\\
 &=& (\B{k}^{-1} \PPsi \W{k} \Sigma_0^{1/2} \eeps)^{\T}
    \B{k}
    (\B{k}^{-1} \PPsi \W{k} \Sigma_0^{1/2} \eeps)\\
 &=& \eeps^{\T} \S \eeps,
\end{eqnarray*}
where the symmetric matrix \( \S \) is defined by \eqref{def S}. Then by the Schur theorem there exist an orthogonal matrix \( \M \) and the diagonal matrix \( \Lam \) composed of the eigenvalues of \( \S \) such that \(     \S = \M^{\T} \Lam \M \). For \( \eeps \sim \norm{0}{I_n} \) and an orthogonal matrix \( \M \) it holds that \( \bar{\eeps} \eqdef \M\eeps \sim \norm{0}{I_n} \). Indeed, \( \EE\M \eeps = \EE \eeps =0 \) and
\begin{equation*}
     \Var \M\eeps = \EE \M\eeps (\M\eeps)^{\T} =
     \M \EE(\eeps \eeps^{\T}) \M^{\T} = \M \M^{\T}  = I_n.
\end{equation*}
Therefore,
\begin{equation*}
    2 \LL ( \W{k}, \mmle{k}, \bbpf{k} )
   \eqdistr \bar{\eeps}^{\T} \Lam \bar{\eeps}\; ,\;\;\;\;
   \bar{\eeps}  \sim \norm{0}{I_n}.
\end{equation*}
On the other hand, the matrix \( \S \) can be written as \( \S = \Sigma_0^{1/2}\W{k}^{1/2}\PPi_k \W{k}^{1/2}\Sigma_0^{1/2} \) with \( \PPi_k = \W{k}^{1/2} \PPsi^{\T} \B{k}^{-1} \PPsi  \W{k}^{1/2}\).
Since \( \PPi_k \) is symmetric and idempotent, i.e. \( \PPi_k^2 = \PPi_k \), it is an orthogonal projector on the linear subspace of dimension \( p = \rank (\B{k}) \) spanned by the rows of \( \PPsi \). Moreover, \( \rank(\PPi_k) = \tr(\PPi_k) = \tr(\W{k}^{1/2} \PPsi^{\T} \B{k}^{-1} \PPsi  \W{k}^{1/2}) = \tr(\B{k}^{-1} \PPsi  \W{k} \PPsi^{\T}) =\tr(\B{k}^{-1} \B{k})= \tr(I_p) = p \). Therefore \( \PPi_k \) has only \( p \) unit eigenvalues and \( n-p \) zero ones. Notice also that the \( n \times n \) matrix \( \S \) has \( \rank(\S) = \rank(\PPi_k \W{k}^{1/2}\Sigma_0^{1/2} ) = \rank(\PPi_k)=p \) as well. Thus  \( 2 \LL ( \W{k}, \mmle{k}, \bbpf{k} ) \eqdistr  \lambda_{1}(\S)\bar\varepsilon_{1}^{2} + \cdots +
    \lambda_{p}(\S)\bar\varepsilon_{p}^{2}\), where  \( \lambda_{1}(\S) , \ldots , \lambda_{p}(\S) \)
are the non-zero eigenvalues of the matrix \( \S  \).

\par Recall the definition of the matrix norm induced by the \( L_2  \) vector norm:
\begin{equation}\label{def om matrix norm}
    \| A \|_{2,in} \eqdef \sqrt{\lambda_{max}(A^{\T}A) }.
\end{equation}
Assumption \( (A3) \) allows to bound the induced \( L_2  \)-norm of the matrix \( \S \):
\begin{eqnarray*}
  \| \S \|_{2,in} &=& \| \Sigma_0^{1/2}\W{k}^{1/2}\PPi_k \W{k}^{1/2}\Sigma_0^{1/2} \|_{2,in}\\
     &\le& \| \Sigma_0^{1/2}\W{k}^{1/2} \|_{2,in} \| \PPi_k \|_{2,in} \| \W{k}^{1/2}\Sigma_0^{1/2} \|_{2,in} \\
     &=&  \lambda_{max}(\W{k}\Sigma_0) \lambda_{max}(\PPi_k)  \\
    &=&   \max_{i}\{ \w{k}{i} \frac{\sigma_{0,i}^2}{\sigma_{i}^2}  \} \\
    &\le & (1 + \delta) \max_{i}\{ \w{k}{i} \} \le  1 + \delta.
\end{eqnarray*}
Therefore, the largest eigenvalue of matrix \( \S \) is bounded: \(\lambda_{max}(\S) \le  1 + \delta \).
\begin{equation*}
    \P \left\{ \lambda_{1}(\S)\bar\eps_1^2 + \cdots
        + \lambda_{p}(\S) \bar\eps_p^2  \ge  \z \right\}
    \le  \P \left\{\lambda_{max}(\S) (\bar\eps_1^2 + \cdots +
        \bar\eps_p^2  ) \ge  \z  \right\}
 \end{equation*}
provides the last assertion.

\end{proof}

\begin{corollary} \textbf{Quasi-parametric risk bounds} \label{param. risk bounds}
\par Let the model be given by \eqref{true model} and \( \bbpf{k} = \bbpf{k}(x) \)
 be defined by \eqref{prameter of BPF}. Assume \( (A3) \). Then for any \( \mu <1/( 1 + \delta) \) we have
\begin{eqnarray}
   \EE \exp \{\mu \LL (\W{k}, \mmle{k}, \bbpf{k} ) \}
    &\le&  \left[1-\mu (1 + \delta) \right]^{-p/2}, \label{param. risk bounds exp}\\
  \EE |2 \LL( \W{k}, \mmle{k}, \bbpf{k} )  |^r
  & \le& (1 + \delta)^r C(p,r)\;, \label{param. risk bounds polinom}
\end{eqnarray}
where
\begin{equation}\label{C(p,r)}
     C(p,r) = \EE|\chi^2_p|^r = 2^r \Gamma (r+ p/2 )/\Gamma(p/2).
\end{equation}
\end{corollary}
\begin{proof}
By \eqref{likelihood = weightsum of sqnormal } and independence of \( \bar{\eps}_i \)
\begin{eqnarray*}
    \EE \exp \{ \mu \LL( \W{k}, \mmle{k}, \bbpf{k} ) \}
    &=&  \EE \exp \left\{ \frac{\mu}{2}
             \sum_{i=1}^{p} \lambda_i(\S) \bar{\eps}_i^2 \right \}\\
    &=&  \prod_{i=1}^{p} \EE \exp \left\{ \frac{\mu}{2} \,
        \lambda_i(\S) \bar{\eps}_i^2  \right \}\\
    &=& \prod_{i=1}^{p} \left[ 1 - \mu \,\lambda_i(\S)  \right]^{-1/2}\\
    &\le & \left[ 1 - \mu \,\lambda_{max}(\S)  \right]^{-p/2} \\
    &\le & [ 1- \mu (1+\delta)]^{-p/2}.
\end{eqnarray*}
Let \( \eta \sim \chi^2_p \). Integration by parts yields the second inequality:
\begin{eqnarray*}
  \EE | 2 \LL( \W{k}, \mmle{k}, \bbpf{k} ) |^r
    &=& \int_0^{\infty}
        \P\left\{2\LL( \W{k}, \mmle{k}, \bbpf{k} )
        \ge \z\right\} r \z^{r-1} \dd\z \\
  &\le& r \int_0^{\infty}
        \P\left\{ \eta \ge \z/(1+\delta)\right\} \z^{r-1} \dd \z \\
  &=& (1 + \delta)^r \,\EE|\eta|^r.
\end{eqnarray*}
\end{proof}

\subsection{Proof of the bounds for the critical values}\label{proof of CV}
Denote for any \( l<k \) the variance of difference \( \mmle{k} - \mmle{l} \) by \( V_{lk} \):
\begin{equation}\label{def_Vlk}
     V_{lk} \eqdef \Var (\mmle{k} - \mmle{l}) \succ 0 .
\end{equation}
Then there exists a unique matrix \( V_{lk}^{1/2} \succ 0 \) such that \(( V_{lk}^{1/2})^2 = V_{lk} \).
\begin{lemma}\label{Tlk large div}
\par Assume \( (A1)- (A4) \). If \( \bbpf{1} = \cdots = \bbpf{k} = \tta \) for \( k \le  K \), then for any \( l <  k \) we have
\begin{eqnarray*}
   \P \left\{2 \LL ( \W{l}, \mmle{l}, \mmle{k} ) \ge \z\right\}
        &\le&
            \P \left\{\eta \ge  \z/\lambda_{max}(V_{lk}^{1/2} \B{l} V_{lk}^{1/2})\right\}\\
        &\le&
            \P \left\{\eta \ge  \z/t_0\right\},\\
   \P \left\{2 \LL ( \W{k}, \mmle{k}, \mmle{l} ) \ge  \z\right\}
   &\le&
     \P \left\{\eta \ge  \z/\lambda_{max}(V_{lk}^{1/2} \B{k} V_{lk}^{1/2})\right\}\\
   &\le&
      \P \left\{\eta \ge  \z/t_1\right\},
 \end{eqnarray*}
 where \( t_0=2(1+\delta) (1 + u_0^{-(k-l)}) \), \( t_1 =2(1+\delta) (1 + u^{(k-l)}) \), \( 1 < u_0 \le u \) are the constants from the assumption \( (A4) \) and   \( \eta \) is a \( \chi^2_p \)-distributed random variable.
\end{lemma}

\begin{proof}
Decomposition~\eqref{linearity if quasiMLE} of \( \mmle{k} \) into deterministic \( \bbpf{k} \) and stochastic parts and the assumption of lemma imply
\begin{equation*}
    \mmle{l}-\mmle{k}
        = \B{l}^{-1}\PPsi\W{l}\Sigma_0^{1/2}\eeps-\B{k}^{-1}\PPsi\W{k}\Sigma_0^{1/2}\eeps
        \eqdistr
            V_{lk}^{1/2}\xi,
\end{equation*}
where \( \xi \) is a standard normal vector in \( \R^{\p} \). Thus by Lemma~\ref{Th. Spokoiny_fitted likelihood} for any \( l<k \)
\begin{equation*}
  2\LL(\W{l}, \mmle{l}, \mmle{k}) =  \| \B{l}^{1/2} (\mmle{l}  - \mmle{k}) \|^2 \\
  \eqdistr  \xi^{\T} V_{lk}^{1/2} \B{l} V_{lk}^{1/2} \xi.
\end{equation*}
 By the Schur theorem there exists an orthogonal matrix \( M \) such that
\begin{equation*}
    \xi^{\T} V_{lk}^{1/2} \B{l} V_{lk}^{1/2} \xi
    \eqdistr  \bar\varepsilon^{\T} M^{\T} \Lambda_{lk} M \bar\varepsilon,
\end{equation*}
where \( \bar\varepsilon \) is a standard normal vector,
\begin{equation*}
     \Lambda_{lk} = \diag (\lambda_{1}(V_{lk}^{1/2} \B{l} V_{lk}^{1/2})) ,
 \cdots , \lambda_{\p}(V_{lk}^{1/2} \B{l} V_{lk}^{1/2}))
\end{equation*}
and \( p = \rank (\B{l})\).
 Therefore,
\begin{equation*}
  2 \LL ( \W{l}, \mmle{l}, \mmle{k} ) \eqdistr
    \lambda_{1}(V_{lk}^{1/2} \B{l} V_{lk}^{1/2})\bar\varepsilon_{1}^{2} + \cdots +
    \lambda_{p}(V_{lk}^{1/2} \B{l} V_{lk}^{1/2})\bar\varepsilon_{p}^{2},
 \end{equation*}
 where \( \lambda_{j}(V_{lk}^{1/2} \B{l} V_{lk}^{1/2}) \), \( j=1,\ldots , p \), are the nonzero eigenvalues
 of \( V_{lk}^{1/2} \B{l} V_{lk}^{1/2} \). Similarly,
 \begin{equation*}
    2 \LL ( \W{k}, \mmle{k}, \mmle{l} ) \eqdistr
    \lambda_{1}(V_{lk}^{1/2} \B{k} V_{lk}^{1/2})\bar\varepsilon_{1}^{2} + \cdots +
    \lambda_{p}(V_{lk}^{1/2} \B{k} V_{lk}^{1/2})\bar\varepsilon_{p}^{2}.
 \end{equation*}

 Denoting by \( \eta \) a \( \chi^2_p \)-distributed random variable we get
 \begin{eqnarray*}
   \P \left\{2 \LL ( \W{l}, \mmle{l}, \mmle{k} ) \ge \z\right\} &\le&
    \P \left\{\eta \ge  \z/\lambda_{max}(V_{lk}^{1/2} \B{l} V_{lk}^{1/2})\right\},\\
   \P \left\{2 \LL ( \W{k}, \mmle{k}, \mmle{l} ) \ge  \z\right\} &\le&
    \P \left\{\eta \ge  \z/\lambda_{max}(V_{lk}^{1/2} \B{k} V_{lk}^{1/2})\right\}.
 \end{eqnarray*}
For any square matrices \( A \) and \( B \) we have
\( (A-B)(A^{\T}-B^{\T}) \preceq 2  (AA^{\T}+BB^{\T}) \).
Applying this bound to the variance of the difference of estimators we obtain
\begin{eqnarray*}
  V_{lk} &=&    \left(\B{l}^{-1}\PPsi\W{l}\Sigma_0^{1/2}-\B{k}^{-1}\PPsi\W{k}\Sigma_0^{1/2} \right)
                \left(\B{l}^{-1}\PPsi\W{l}\Sigma_0^{1/2}-\B{k}^{-1}\PPsi\W{k}\Sigma_0^{1/2}\right)^{\T} \\
&\preceq& 2 (\B{l}^{-1}\PPsi\W{l} \Sigma_0 \W{l}  \PPsi^{\T}\B{l}^{-1} + \B{k}^{-1}\PPsi\W{k} \Sigma_0 \W{k}  \PPsi^{\T}\B{k}^{-1} ) \\
&=&  2V_l +2V_k,
\end{eqnarray*}
where \( V_l=\Var\mmle{l} \), \( l\le k \). By \eqref{Vk bound} and Assumption \( (A4) \) we have
 \begin{eqnarray*}
& & V_{l}  \preceq  (1+ \delta)\B{l}^{-1},\\
& & V_{k}  \preceq  (1+ \delta)\B{k}^{-1} \preceq  (1+ \delta)  u_0^{-(k-l)}  \B{l}^{-1}, \\
& & V_{lk}  \preceq 2(1+\delta) (1 + u_0^{-(k-l)}) \B{l}^{-1}.
\end{eqnarray*}
Therefore,
\begin{equation}\label{bound for B_l}
  \B{l} \preceq 2(1+\delta) (1 + u_0^{-(k-l)}) V_{lk}^{-1}.
 \end{equation}
This provides the following bound:
 \begin{eqnarray}
\lambda_{max}(V_{lk}^{1/2} \B{l} V_{lk}^{1/2})
&=& \sup_{\| \gamma \| = 1} \gamma^{\T}  V_{lk}^{1/2} \B{l} V_{lk}^{1/2} \gamma  \nn
&\le & 2(1+\delta) (1 + u_0^{-(k-l)})\label{lambda max l}.
\end{eqnarray}
Similarly,
\begin{eqnarray}
 V_{lk} &\preceq& 2(1+\delta) (1 + u^{(k-l)}) \B{k}^{-1} ,\nn
  \lambda_{max}(V_{lk}^{1/2} \B{k} V_{lk}^{1/2})&\le & 2(1+\delta) (1 + u^{(k-l)})
                                                                    \label{lambda max k}.
\end{eqnarray}
These bounds imply
\begin{eqnarray*}
   \P \left\{2 \LL ( \W{l}, \mmle{l}, \mmle{k} ) \ge \z\right\}
        &\le&
            \P \left\{\eta \ge  \z/\lambda_{max}(V_{lk}^{1/2} \B{l} V_{lk}^{1/2})\right\}\\
        &\le&
            \P \left\{\eta \ge  \z\left[2(1+\delta) (1 + u_0^{-(k-l)})\right]^{-1}\right\}\\
   \P \left\{2 \LL ( \W{k}, \mmle{k}, \mmle{l} ) \ge  \z\right\}
   &\le&
     \P \left\{\eta \ge  \z/\lambda_{max}(V_{lk}^{1/2} \B{k} V_{lk}^{1/2})\right\}\\
   &\le&
      \P \left\{\eta \ge  \z\left[2(1+\delta) (1 + u^{(k-l)})\right]^{-1}\right\}
 \end{eqnarray*}

\end{proof}

\begin{lemma} \label{expmoment corol_new}
Under the conditions of preceding lemma for any \( l<k \), \( \mu_0 < t_0^{-1} \),
\( \mu_1 < t_1^{-1} \) we have
\begin{eqnarray*}\label{expmoment_new}
\EE \exp \{ \mu_0 \LL(\W{l}, \mmle{l}, \mmle{k}) \}
        &\le&
        [1 - \mu_0 t_0]^{-p/2},\\
        \EE \exp \{ \mu_1 \LL(\W{k}, \mmle{k}, \mmle{l}) \}
        &\le &
        [1 - \mu_1 t_1]^{-p/2},
\end{eqnarray*}
where \( t_0=2(1+\delta) (1 + u_0^{-(k-l)}) \), \( t_1 =2(1+\delta) (1 + u^{(k-l)}) \) and the constants \( 1 < u_0 \le u \) are from Assumption \( (A4) \).
\end{lemma}
\begin{proof}
The  statement of the lemma is justified similarly to the proof of Corollary~\ref{param. risk bounds}. The bounds \eqref{lambda max l} and \eqref{lambda max k} imply the bounds for the corresponding moment generating functions:
\begin{eqnarray*}
\EE \exp \{ \mu \LL(\W{l}, \mmle{l}, \mmle{k}) \}
        &=& \prod_{j=1}^{p} \EE \exp \{ \frac{\mu}{2} \lambda_{j}(V_{lk}^{1/2} \B{l} V_{lk}^{1/2})
            \bar\varepsilon_{j}^2\}\\
        &=&    \prod_{j=1}^{p} [1- \mu  \lambda_{j}(V_{lk}^{1/2} \B{l} V_{lk}^{1/2}) ]^{-1/2}\\
        &\le&
        [1 - \mu \lambda_{max}(V_{lk}^{1/2} \B{l} V_{lk}^{1/2}) ]^{-p/2}\\
          &\le&
        [1 - 2\mu (1+\delta) (1 + u_0^{-(k-l)})]^{-p/2}\;,\\
\EE \exp \{ \mu \LL(\W{k}, \mmle{k}, \mmle{l}) \}
        &\le &
        [1 - \mu \lambda_{max}(V_{lk}^{1/2} \B{k} V_{lk}^{1/2})]^{-p/2}\\
         &\le & [1 - 2\mu (1+\delta) (1 + u^{(k-l)})]^{-p/2}.
\end{eqnarray*}
\end{proof}

\begin{lemma}\label{polmoment corol_new}
Under the conditions of Lemma~\ref{Tlk large div} for any \( l<k \) we have
\begin{eqnarray*}
    \EE| 2\LL(\W{l}, \mmle{l}, \mmle{k}) |^r &\le& 2^r C(p,r) (1+\delta)^r (1 + u_0^{-(k-l)})^r,\\
    \EE| 2\LL(\W{k}, \mmle{k}, \mmle{l}) |^r &\le& 2^r C(p,r) (1+\delta)^r (1 + u^{(k-l)})^r,
\end{eqnarray*}
where \(  C(p,r) \) is given by \eqref{C(p,r)}.
\end{lemma}
\begin{remark}
The RHS's of Lemmas~\ref{expmoment corol_new} and \ref{polmoment corol_new} are highly asymmetric. Recall that here \( \bbpf{1} = \cdots = \bbpf{k} = \tta \), \( l<k \) and \( 1 < u_0 \le u \). The bounds for the log-likelihood ratio corresponding to the \( l \)-th scale \( \LL(\W{l}, \mmle{l}, \mmle{k}) \) are close to the bounds for their parametric counterpart \( \LL(\W{l}, \mmle{l}, \tta) \) given by Corollary~\ref{param. risk bounds}. It is not surprising because, if the parametric model is satisfied up to the scale \( k \), for the MLE \( \mmle{k} \) more data were used and the estimator \( \mmle{k} \) w.r.t. \( \mmle{l} \) acts approximately as the true parameter \( \tta \). On the contrary, the risk bounds for \( \LL(\W{k}, \mmle{k}, \mmle{l}) \) are quite large since for the larger \( k \)-th scale \( \mmle{l} \) is a bad estimator with large variance.
\end{remark}
\begin{proof} Integration by parts and Lemma \ref{Tlk large div} yield for the second assertion
\begin{eqnarray*}
  \EE| 2\LL(\W{k}, \mmle{k}, \mmle{l}) |^r
    &=&
    r \int_{0}^{\infty}  \P \left\{2 \LL ( \W{k}, \mmle{k}, \mmle{l} ) \ge   \z\right\} \z^{r-1} \dd \z \\
   &\le &
   r \int_{0}^{\infty}    \P \left\{\eta \ge  \z \left[2(1+\delta)
        (1 + u^{(k-l)})\right]^{-1}  \right\} \z^{r-1} \dd \z \\
   &=& 2^r  (1+\delta)^r (1 + u^{(k-l)})^r \E|\eta|^r,
\end{eqnarray*}
where \( \eta \sim  \chi^2_p\). The first assertion is proved similarly.
\end{proof}

\begin{proof} \emph{of Theorem \ref{upper bound} Theoretical choice of the critical values.}
The risk corresponding to the adaptive estimate can be represented as a sum of risks of
the false alarms at each step of the procedure:
\begin{equation}\label{proof_CV_decomposition}
\EE_{0, \Sigma} |(\mmle{k} - \aadapest_k )^{\T} \B{k} (\mmle{k} - \aadapest_k ) |^r
    =
    \sum_{m=1}^{k-1} \EE_{0, \Sigma} |(\mmle{k} - \mmle{m} )^{\T} \B{k} (\mmle{k} - \mmle{m} ) |^r \I{\{ \aadapest_{k} = \mmle{m}\}}.
\end{equation}

By the definition of the last accepted estimate \( \aadapest_{k} \), for any \( m =1, \ldots, k-1 \), the event \( \{ \aadapest_{k} = \mmle{m} \} \) happens if for some \( l = 1, \ldots ,m \) the statistic \( T_{l, m+1} >\z_l \). Thus
\begin{equation*}
    \{ \aadapest_{k} = \mmle{m} \} \subseteq \bigcup_{l=1}^{m} \{ T_{l, m+1} > \z_l \}.
\end{equation*}
It holds also that for any positive \( \mu \)
\begin{eqnarray*}
    \I{\{ T_{l, m+1} > \z_l \}}
        &=&
         \I{\{ 2 \LL(\W{l}, \mmle{l}, \mmle{m+1}) - \z_l > 0\}}\\
        &\le & \exp \{ \frac{\mu}{2} \LL(\W{l}, \mmle{l}, \mmle{m+1}) - \frac{\mu}{4} \z_l\}.
\end{eqnarray*}
This and the Cauchy-Schwarz  inequality imply for \( m = 1, \ldots, k-1 \) the following bound:
\begin{eqnarray}\label{proof_CV_decomposition_CauchySchwarz}
 && \EE_{0, \Sigma} |(\mmle{k} - \mmle{m} )^{\T} \B{k} (\mmle{k} - \mmle{m} ) |^r \I{\{ \aadapest_{k} = \mmle{m} \}}\\
   &=& \EE_{0, \Sigma} |2 \LL (\W{k}, \mmle{k}, \mmle{m})|^r \I{\{ \aadapest_{k} = \mmle{m}\}} \nn
  &\le & \sum_{l=1}^{m} e^{-\frac{\mu}{4}\z_l} \EE_{0, \Sigma} \left[|2 \LL (\W{k}, \mmle{k}, \mmle{m})|^r
    \exp{\{ \frac{\mu}{2} \LL(\W{l}, \mmle{l}, \mmle{m+1}) \}}\right] \nn \nonumber
    &\le & \sum_{l=1}^{m} e^{-\frac{\mu}{4}\z_l}
        \left\{  \EE_{0, \Sigma} \left[ |2 \LL (\W{k}, \mmle{k}, \mmle{m})|^{2r} \right]
          \right\}^{\frac{1}{2} }
    \left\{ \EE_{0, \Sigma} \left[ \exp{\{ \mu \LL(\W{l}, \mmle{l}, \mmle{m+1}) \}} \right]\right\}^{\frac{1}{2}}.
\end{eqnarray}
By the first statement of Lemma \ref{expmoment corol_new} with \( \delta = 0 \)
\begin{equation*}
    \EE_{0, \Sigma} \left[ \exp{\{ \mu \LL(\W{l}, \mmle{l}, \mmle{m+1}) \}} \right]
        \le [ 1- 2\mu (1+u_0^{-(m+1-l)})]^{-\frac{p}{2}}
\end{equation*}
for any \( \mu < [2 (1+u_0^{-(m+1-l)}]^{-1} \). Since \( u_0>1 \) we have \( [2 (1+u_0^{-(m+1-l)}]^{-1} >1/4 \) and the statement is valid for any \( \mu \in (0,1/4) \). Inequality \( [ 1- 2\mu (1+u_0^{-(m+1-l)})]^{-p/2} < [1-4 \mu]^{-p/2} \) provides for any \( \mu \in (0,1/4) \)
\begin{equation}\label{proof_CV_bound for exp.mom.}
    \EE_{0, \Sigma} \left[ \exp{\{ \mu \LL(\W{l}, \mmle{l}, \mmle{m+1}) \}} \right]
        < ( 1- 4\mu)^{-p/2}.
\end{equation}
By the second statement of Lemma \ref{polmoment corol_new}
\begin{equation}\label{proof_CV_bound for pol.mom.}
    \EE_{0, \Sigma}  |2 \LL (\W{k}, \mmle{k}, \mmle{m})|^{2r} \le C(p,2r) 2^{2r} (1+u^{k-m})^{2r}.
\end{equation}

Putting together \eqref{proof_CV_decomposition}, \eqref{proof_CV_decomposition_CauchySchwarz}, \eqref{proof_CV_bound for exp.mom.} and \eqref{proof_CV_bound for pol.mom.}  we obtain
\begin{eqnarray*}
   && \EE_{0, \Sigma} |(\mmle{k} - \aadapest_k )^{\T} \B{k} (\mmle{k} - \aadapest_k ) |^r \\
  &\le & 2^r \sqrt{C(p,2r)}   ( 1- 4\mu)^{-p/4}
   \sum_{m=1}^{k-1} \sum_{l=1}^{m} e^{-\frac{\mu}{4}\z_l}(1+u^{k-m})^r\\
   &=& 2^r \sqrt{C(p,2r)}  ( 1- 4\mu)^{-p/4}
        \sum_{l=1}^{k-1} e^{-\frac{\mu}{4}\z_l} \sum_{m=l}^{k-1}(1+u^{k-m})^r\\
   &\le & 2^{2r} \sqrt{C(p,2r)}   ( 1- 4\mu)^{-p/4} (1-u^{-r})^{-1}
        \sum_{l=1}^{k-1} e^{-\frac{\mu}{4}\z_l} u^{r(k-l)},
\end{eqnarray*}
because \( -(k-l) < -(m-l) \) and
\begin{eqnarray*}
    \sum_{m=l}^{k-1}(1+u^{(k-m)})^r  &=& u^{r(k-l)} \sum_{m=l}^{k-1}(u^{-(k-l)}+u^{-(m-l)})^r\\
    &<& 2^r u^{r(k-l)} \sum_{m=l}^{k-1} u^{-r(m-l)}\\
    &<& 2^r u^{r(k-l)} (1-u^{-r})^{-1}.
    \end{eqnarray*}
Since \( u^{r(k-l)} \le u^{r(K-l)} \) for any \( l<k\le K \) the choice of the threshold of the form
\begin{equation*}
    \z_l = \frac{4}{\mu} \left\{ r(K-l) \log u + \log{(K/\alpha)} - \frac{p}{4} \log(1-4\mu )
    - \log(1-u^{-r}) + \bar{C}(p,r) \right\}
\end{equation*}
with an arbitrary constant \( \mu \in (0,1/4) \), \( u>1 \) from Assumption \( (A4) \), \( r>0 \) and \( \alpha \in (0,1] \) from the PC's and with
\begin{equation*}
    \bar{C}(p,r) =\log\left\{ \frac{2^{2r} [\Gamma(2r + p/2) \Gamma(p/2)]^{1/2}}{\Gamma(r + p/2)} \right \}
\end{equation*}
provides the required by PC bounds
\begin{equation*}
    \EE_{0,\Sigma} |(\mmle{l} - \aadapest_{l})^{\T}
        \B{l} (\mmle{l} - \aadapest_{l})|^{r}
        \leq \alpha  C(p,r)\; \; \;
        \text{for all}\;\; l=2, \ldots, K.
 \end{equation*}
\end{proof}

\subsection{Matrix results}\label{Matrixresults}
\begin{lemma}\label{simidefinitness of SSigma }
The matrices \(J_k \otimes \Sigma\) and \(J_k \otimes \Sigma_0\) are positive semidefinite for any \( k=2, \ldots, K \).
\par Moreover, under Assumption \( (A3)\) with the same \( \delta \), the similar to
\( (A3)\) relation holds for the covariance matrices \( \SSigma_{k}   \) and \( \SSigma_{k,0} \) of linear estimates:
\begin{equation*}
  (1-\delta) \SSigma_{k} \preceq  \SSigma_{k,0} \preceq (1+\delta)  \SSigma_{k}\;,\;\; k \le K.
\end{equation*}
\end{lemma}
\begin{proof} 
\par
Symmetry of \( J_k\) and \(\Sigma \), (respectively, \(\Sigma_0 \) ) implies symmetry of \(J_k \otimes \Sigma\), (respectively, \(J_k \otimes \Sigma_0\)). Notice that any vector \( \gamma_{nk} \in \R^{nk} \) can be represented as a partitioned vector \( \gamma_{nk}^{\T} = ((\gamma_{nk}^{(1)})^{\T},(\gamma_{nk}^{(2)})^{\T}, \ldots, (\gamma_{nk}^{(k)})^{\T}) \), with \( \gamma_{nk}^{(l)} \in \R^{n} \), \( l=1, \ldots, k \). Then
\begin{equation}\label{gamma nk}
    \gamma_{nk}^{\T} (J_k \otimes \Sigma) \gamma_{nk}
        = \big(\sum_{l=1}^k \gamma_{nk}^{(l)}\big)^{\T}  \Sigma \big( \sum_{l=1}^k \gamma_{nk}^{(l)} \big)
        = \tilde{\gamma}_n^{\T} \, \Sigma \, \tilde{\gamma}_n,
\end{equation}
where \( \tilde{\gamma}_n \eqdef \sum_{l=1}^k \gamma_{nk}^{(l)} \in \R^n \). Because \( \Sigma \succ 0 \) it implies \( \tilde{\gamma}_n^{\T} \Sigma \,\tilde{\gamma}_n > 0 \) for all \( \tilde{\gamma}_n \ne 0 \). But even for \( \gamma_{nk} \ne 0 \), if its subvectors \( \{ \gamma_{nl}^{(l)} \} \) are linearly dependent, \( \tilde{\gamma}_n \) can be zero. Thus there exists a nonzero vector \( \gamma \) such that \( \gamma^{\T} (J_k \otimes \Sigma) \gamma =0 \). This means positive semidefiniteness.

\par The second assertion follows from the observation that
Assumption \( (A3)\) due to the equality \eqref{gamma nk} also holds for the Kronecker product
\begin{equation}\label{ineq for the midle parts of SSigmas}
 (1-\delta) J_k \otimes \Sigma \preceq J_k \otimes \Sigma_0 \preceq (1+\delta) J_k \otimes \Sigma.
\end{equation}
Therefore
\begin{equation*}
     (1-\delta) \DD_k(J_k \otimes \Sigma) \DD_k^{\T}
        \preceq \DD_k(J_k \otimes \Sigma_0 ) \DD_k^{\T}
        \preceq (1+\delta) \DD_k(J_k \otimes \Sigma) \DD_k^{\T}.
\end{equation*}

\end{proof}

\begin{lemma}\label{nonsingularity of SSigma rectang}
Fix \( x \in \R^d \). Suppose that the weights \( \{ \w{l}{\,i}(x) \} \) satisfy
\begin{equation}\label{binar product}
    \w{l}{\,i}(x) \w{m}{\,i}(x) = \w{l}{\,i}(x) \;,\;\; l \le m.
\end{equation}
Then under Assumptions \( (A1) \), \( (A2) \), \( (A4) \) the covariance matrix \( \SSigma_k \) defined by \eqref{def SSigma} is nonsingular with
\begin{equation}\label{det of SSigma rect}
    \det \SSigma_k = \det \B{k}^{-1} \prod_{l=2}^k \det (\B{l-1}^{-1} -\B{l}^{-1}) > 0
                    \;,\;\;k=2, \ldots, K.
\end{equation}
\end{lemma}

\begin{remark}
The condition \eqref{binar product} holds for rectangular kernels with nested supports.
\end{remark}

\begin{proof}
The condition \eqref{binar product} implies
  \begin{equation*}
     \W{l} \Sigma \W{m}
     = \diag(\w{l}{1} \w{m}{1}/\sigma_1^2 , \ldots, \w{l}{n} \w{m}{n}/\sigma_n^2 ) = \W{l}
\end{equation*}
for any \( l \le m \). Thus the blocks of \( \SSigma_k \) simplify to
\begin{equation*}
     D_l \Sigma D_m^{\T} =
\B{l}^{-1}  \PPsi  \W{l} \Sigma \W{m}  \PPsi^{\T}  \B{m}^{-1} =
\B{l}^{-1}  \PPsi  \W{l}   \PPsi^{\T}  \B{m}^{-1}
\end{equation*}
and \( \SSigma_k  \) has a simple structure:
\[
\SSigma_k =
\begin{pmatrix}
    \B{1}^{-1} & \B{2}^{-1}  & \B{3}^{-1}  & \ldots &       \B{k}^{-1} \\
    \B{2}^{-1} & \B{2}^{-1}  & \B{3}^{-1}  & \ldots &       \B{k}^{-1} \\
    \vdots     &   \vdots    &   \vdots    &  \vdots &       \vdots   \\
    \B{k}^{-1} & \B{k}^{-1}  & \B{k}^{-1}  & \ldots &        \B{k}^{-1}
\end{pmatrix}.
\]
Then the determinant of \( \SSigma_k \) coincides with the determinant of the following irreducible block triangular matrix:
\[
\det \SSigma_k =
\begin{vmatrix}
    \B{1}^{-1}-\B{2}^{-1} &  \B{2}^{-1}-\B{3}^{-1}  & \ldots & \B{k-1}^{-1}-\B{k}^{-1} & \B{k}^{-1} \\
     \bb{0}               &  \B{2}^{-1}-\B{3}^{-1}  & \ldots & \B{k-1}^{-1}-\B{k}^{-1} & \B{k}^{-1} \\
     \vdots               &          \vdots         & \vdots &       \vdots            &   \vdots   \\
     \bb{0}               &   \bb{0}                & \ldots & \B{k-1}^{-1}-\B{k}^{-1} & \B{k}^{-1} \\
     \bb{0}               &   \bb{0}                & \bb{0} &      \bb{0}             & \B{k}^{-1}
\end{vmatrix}
\] implying
\begin{equation*}
    \det \SSigma_k = \det(\B{1}^{-1}-\B{2}^{-1}) \det(\B{2}^{-1}-\B{3}^{-1})\cdot \ldots \cdot
        \det(\B{k-1}^{-1}-\B{k}^{-1}) \det\B{k}^{-1}.
\end{equation*}
Clearly the matrix \( \SSigma_k \) is nonsingular if all the matrices \( \B{l-1}^{-1}-\B{l}^{-1} \) are nonsingular. By \( (A1) \) and \( (A2) \) \( \B{l} \succ 0 \) for any \( l \). By \( (A4) \) there exists \( u_0> 1 \) such that \( \B{l}  \succeq  u_0  \B{l-1} \), therefore
\( \B{l-1}^{-1}-\B{l}^{-1} \succeq (1-1/u_0) \B{l-1}^{-1} \succ \B{l-1}^{-1} \succ 0\).
\end{proof}

\begin{lemma}\label{MGF for joint distribution}
In the ``nonparametric situation'' the moment generation function (mgf) of the joint distribution of \( \mmle{1}, \ldots , \mmle{K} \) is
\begin{equation}\label{MGF under alternative}
    \EE \exp \big\{ \gamma^{\T} ( \vec \tilde \TTa_K - \vec \TTa^*_K) \big\}
        = \exp \bigg\{ \frac{1}{2} \gamma^{\T} \SSigma_{K,0} \, \gamma \bigg\}.
\end{equation}
Thus, provided that \( \SSigma_{K,0} \succ 0 \), it holds that \( \vec \tilde \TTa_K \sim \norm{\vec \TTa^*_K}{\SSigma_{K,0}} \).
\par Similarly, in the ``parametric situation'', if \( \SSigma_K \succ 0 \), then the joint distribution of \( \vec \tilde \TTa_K \) is \(  \norm{\vec \TTa_K}{\SSigma_K} \) with the mgf:
\begin{equation}\label{MGF under null}
    \EE \exp \big\{ \gamma^{\T} ( \vec \tilde \TTa_K - \vec \TTa_K) \big\}
        = \exp \bigg\{ \frac{1}{2} \gamma^{\T} \SSigma_K \, \gamma \bigg\}.
\end{equation}
\end{lemma}
\begin{proof}
Let \( \gamma \in \R^{pK} \) be written in a partitioned form
\( \gamma^{\T} = (\gamma_1^{\T}, \ldots, \gamma_K^{\T}) \)
with \( \gamma_l \in \R^p \), \( l=1, \ldots, K \). Then the mgf for the centered random vector
\( \vec \tilde \TTa_K - \vec \TTa^*_K \in \R^{pK} \), due to the decomposition \eqref{linearity if quasiMLE} \( \mmle{l} = \bbpf{l} + D_l \Sigma_0^{1/2} \eeps\) with \( D_l = \B{l}^{-1}\PPsi \W{l}\),
can be represented as follows:
  \begin{eqnarray*}
    &&\EE \exp \big\{  \gamma^{\T} ( \vec \tilde \TTa_K - \vec \TTa^*_K) \big\}
        = \EE \exp\big\{ \sum_{l=1}^K \gamma_l^{\T} (\mmle{l} - \bbpf{l}) \big\} \\
    &=& \EE \exp\big\{ \sum_{l=1}^K \gamma_l^{\T} D_{l} \Sigma_0^{1/2} \eeps \big\}
        = \EE \exp\big\{ \big(\sum_{l=1}^K D_{l}^{\T} \gamma_l  \big)^{\T} \Sigma_0^{1/2} \eeps \big\}.
  \end{eqnarray*}
A trivial observation that \( \sum_{l=1}^K D_{l}^{\T} \gamma_l   \) is a vector in \( \R^n \) and
\( \Sigma_0^{1/2} \eeps \sim \norm{0}{\Sigma_0} \) by \eqref{true model} implies by the definition of \( \SSigma_{K,0} \) the first assertion of the lemma, because
\begin{eqnarray*}
  &&\EE \exp\big\{ \big(\sum_{l=1}^K D_{l}^{\T} \gamma_l  \big)^{\T} \Sigma_0^{1/2} \eeps \big\}
     = \exp \bigg\{ \frac{1}{2} \big(\sum_{l=1}^K D_{l}^{\T} \gamma_l  \big)^{\T}
            \Sigma_0 \big(\sum_{l=1}^K D_{l}^{\T} \gamma_l  \big) \bigg\}  \\
  &=&\exp \bigg\{ \frac{1}{2}\big(\DD_K^{\T} \gamma \big)^{\T} (J_K \otimes \Sigma_0) \DD_K^{\T} \gamma \bigg\} = \exp \bigg\{ \frac{1}{2} \gamma^{\T} \SSigma_{K,0} \, \gamma \bigg\},
\end{eqnarray*}
here \( \DD_K \) is defined by \eqref{def DD_K}.
\end{proof}

\subsection{Proof of the propagation property} \label{proof propagation}

\begin{lemma}\label{KL for joint distr}
The Kullback-Leibler divergence between the distributions of \( \vec \tilde \TTa_k \) under the true measure and under the ``parametric'' has the following form:
\begin{eqnarray}
  & &2\KL(\P_{\ff, \Sigma_0}^k,\P_{\tta, \Sigma}^k) \eqdef 2\EE_{\ff, \Sigma_0}
  \log\big( \frac{\dd \P_{\ff, \Sigma_0}^k}{\dd \P_{\tta, \Sigma}^k}  \big) \label{KL joint}\nn
  &=&  \Delta(k)  +\log\bigg(\frac{\det \SSigma_k}{\det \SSigma_{k,0}} \bigg)
            + \tr (\SSigma_k^{-1} \SSigma_{k,0}) -pk,
\end{eqnarray}
where
\begin{eqnarray}
  b(k) &\eqdef& \vec \TTa^*_k - \vec \TTa_k \label{def_b(k)} \\
  \Delta(k) & \eqdef & b(k)^{\T} \SSigma_k^{-1} b(k)  \label{def_Delta(k)}.
\end{eqnarray}

\end{lemma}
\begin{proof} 

\par Denote the Radon-Nikodym derivative by \( Z_k \eqdef \dd \P_{\ff, \Sigma_0}^k / \dd \P_{\tta, \Sigma}^k \). Then
\begin{eqnarray}\label{logdp}
    \log\big(Z_k(y)\big) =
        \frac{1}{2} \log\bigg(\frac{\det \SSigma_k}{\det \SSigma_{k,0}} \bigg)
        &-&\frac{1}{2} \| \SSigma_{k,0}^{-1/2} (y - \vec \TTa^*_k) \|^2  \nn
        &+& \frac{1}{2}\| \SSigma_k^{-1/2} (y - \vec \TTa_k) \|^2
\end{eqnarray}
 can be considered as a quadratic function of \( \vec \TTa_k \). By the Taylor expansion at the point \( \vec \TTa^*_k \) the last expression reads as follows
\begin{eqnarray*}
 & & \log\big(Z_k(y)\big) =
        \frac{1}{2} \log\bigg(\frac{\det \SSigma_k}{\det \SSigma_{k,0}} \bigg)
        -\frac{1}{2} \| \SSigma_{k,0}^{-1/2} (y - \vec \TTa^*_k) \|^2  \nn
&+&  \frac{1}{2} \| \SSigma_k^{-1/2} (y - \vec \TTa^*_k) \|^2
+ b(k)^{\T} \SSigma_k^{-1} (y- \vec \TTa^*_k) + \frac{1}{2} \Delta(k).
\end{eqnarray*}
Then the expression for the Kullback-Leibler divergence can be written in the following way:
\begin{eqnarray*}
  & &\KL(\P_{\ff, \Sigma_{0}}^k,\P_{\tta, \Sigma}^k) \eqdef \EE_{\ff, \Sigma_0}
  \log\big(Z_k\big) \\
  &=&  \frac{1}{2} \log\bigg(\frac{\det \SSigma_k}{\det \SSigma_{k,0}} \bigg)
       + \frac{1}{2} \Delta(k)
       + \frac{1}{2} \EE \big\{ \| \SSigma_k^{-1/2} \SSigma_{k,0}^{1/2} \xi \|^2 -\| \xi \|^2
       + 2 b(k)^{\T} \SSigma_k^{-1} \SSigma_{k,0}^{1/2} \xi  \big\},
\end{eqnarray*}
where \( \xi \sim \norm{0}{I_{pk}} \). This implies
\begin{equation}\label{KL general}
    2 \KL(\P_{\ff, \Sigma_0}^k,\P_{\tta, \Sigma}^k)
        = \Delta(k)
            +\log\bigg(\frac{\det \SSigma_k}{\det \SSigma_{k,0}} \bigg)
            + \tr (\SSigma_k^{-1} \SSigma_{k,0}) -pk.
\end{equation}
\par In the case of \emph{homogeneous errors} with \( \sigma_{0,i} = \sigma_{0} \) and
 \( \sigma_{i} = \sigma, i = 1, \ldots ,n \) the calculations simplify a lot. Now
 \begin{equation*}
    \SSigma_k = \sigma^2 \VV_k, \;\;\; \SSigma_{k,0} = \sigma^2_0 \VV_k
 \end{equation*}
with a \( pk \times pk  \) matrix \( \VV_k \) defined as
\begin{equation*}
     \VV_k = \big( \bar{D}_1 \oplus  \cdots \oplus \bar{D}_k  \big)
                    \big( J_k \otimes I_n \big)
                    \big( \bar{D}_1 \oplus  \cdots \oplus \bar{D}_k  \big)^{\T},
\end{equation*}
where \( \bar{D}_l = (\PPsi \cc{W}_l \PPsi^{\T})^{-1} \PPsi \cc{W}_l  \), \( l= 1, \ldots, k \) does not depend on \( \sigma \). Then \( \Delta(k) = \sigma^{-2} \Delta_1(k)  \), with \( \Delta_1(k) \eqdef b(k)^{\T} \VV_k^{-1} b(k) \), \( \det \SSigma_k / \det\SSigma_{k,0}  = (\sigma^2/\sigma_0^2)^{pk} \), and the expression for the Kullback-Leibler divergence reads as follows:
\begin{eqnarray}\label{KL hom}
    \KL(\P_{\ff, \Sigma_{0}}^k,\P_{\tta, \Sigma}^k)
        &=& p k \log\big( \frac{\sigma}{\sigma_0} \big)
            + \frac{1}{2} \Delta(k) + \frac{p k}{2} \big( \frac{\sigma_0^2}{\sigma^2} - 1 \big)\\
            \nonumber
        &=&  p k \log\big( \frac{\sigma}{\sigma_0} \big)
            + \frac{1}{2 \sigma^2} b(k)^{\T} \VV_k^{-1} b(k)
            + \frac{p k}{2} \big( \frac{\sigma_0^2}{\sigma^2} - 1 \big) ,
\end{eqnarray}
implying the same asymptotic behavior as in \eqref{bounds for KL}.

\end{proof}

\begin{proof} \emph{of Theorem \ref{Propagation result theorem} (Propagation property)}
\par Notice that for any nonnegative measurable function \( g = g(\tilde \TTa_k) \) the Cauchy-Schwarz inequality implies
\begin{equation}\label{changing measure}
    \EE_{\ff, \Sigma_0} [g] = \EE_{\tta, \Sigma} [g Z_k]
    \le \big(\EE_{\tta, \Sigma} [g^2] \big)^{1/2} \big(\EE_{\tta, \Sigma} [Z_k^2] \big)^{1/2}
\end{equation}
with the Radon-Nikodym derivative \(     Z_k = \dd \P^k_{\ff, \Sigma_0}/\dd \P^k_{\tta, \Sigma}\). One gets the first assertion taking \( g=|(\mmle{k} - \tta)^{\T} \B{k} (\mmle{k} - \tta)|^{r/2} \), and applying ``the parametric risk bound'' with \( \delta=0 \) from \eqref{param. risk bounds polinom}:
\begin{eqnarray*}
  \EE[g] &\le& \big(\EE_{\tta, \Sigma} |(\mmle{k} - \tta)^{\T} \B{k} (\mmle{k} - \tta)|^{r} \big)^{1/2} \big(\EE_{\tta, \Sigma} [Z_k^2] \big)^{1/2} \\
   &=& \big(\EE_{\tta, \Sigma} |2\LL(\W{k}, \mmle{k}, \tta )|^r \big)^{1/2} \big(\EE_{\tta, \Sigma} [Z_k^2] \big)^{1/2}\\
   &\le&  (\EE|\chi^2_p|^r)^{1/2} \big(\EE_{\tta, \Sigma} [Z_k^2] \big)^{1/2}.
\end{eqnarray*}
The second assertion of the theorem is treated similarly by application of the pivotality property from Lemma~\ref{Pivotality property} and the propagation conditions \eqref{PC}.

To calculate \( \EE_{\tta, \Sigma} [Z_k^2] \) let us consider \( \log Z_k \) given by
\begin{eqnarray*}
    \log\big(Z_k(y)\big) =
        \frac{1}{2} \log\bigg(\frac{\det \SSigma_k}{\det \SSigma_{k,0}} \bigg)
        &-&\frac{1}{2} \| \SSigma_{k,0}^{-1/2} (y - \vec \TTa^*_k) \|^2  \nn
        &+& \frac{1}{2}\| \SSigma_k^{-1/2} (y - \vec \TTa_k) \|^2
\end{eqnarray*}
as a function of \( \vec \TTa^*_k \). Application of the Taylor expansion at the point \(  \vec \TTa_k \) yields
\begin{eqnarray*}
  2\log Z_k &=& \log \frac{\det \SSigma_k}{\det \SSigma_{k,0}}
                 - \| \SSigma_{k,0}^{-1/2} (y - \vec \TTa_k) \|^2
                 + \| \SSigma_k^{-1/2} (y - \vec \TTa_k) \|^2 \\
            &+&  2 b(k)^{\T} \SSigma_{k,0}^{-1}(y - \vec \TTa_k)
                 - b(k)^{\T} \SSigma_{k,0}^{-1} b(k).
\end{eqnarray*}
With \( \xi \sim \norm{0}{I_{pk}} \) the second moment of the Radon-Nikodym derivative reads as follows
\begin{eqnarray}
   & & \EE_{\tta, \Sigma} [Z_k^2] \nn
   &=&     \frac{\det \SSigma_k}{\det \SSigma_{k,0}}
            \exp\{ - b(k)^{\T} \SSigma_{k,0}^{-1} b(k)\}
            \EE \exp \{ -\| \SSigma_{k,0}^{-1/2} \SSigma_k^{1/2} \xi \|^2
                        + \| \xi \|^2
                        + 2 b(k)^{\T} \SSigma_{k,0}^{-1} \SSigma_k^{1/2} \xi \}\nn
   &=&    \frac{\det \SSigma_k}{\det \SSigma_{k,0}}
                \big[\det \big(2 \SSigma_k^{1/2}\SSigma_{k,0}^{-1}\SSigma_k^{1/2} - I_{pk} \big )\big]^{-1/2}\nn
   &\times&       \exp\{ 2 b(k)^{\T} \SSigma_{k,0}^{-1} \SSigma_k^{1/2} \big(2 \SSigma_k^{1/2} \SSigma_{k,0}^{-1}  \SSigma_k^{1/2} -I_{pk}\big)^{-1} \SSigma_k^{1/2} \SSigma_{k,0}^{-1} b(k)
           - b(k)^{\T} \SSigma_{k,0}^{-1} b(k)\} \nn
   &=&    \frac{\det \SSigma_k}{\det \SSigma_{k,0}}
                \big[ \prod_{j=1}^{pk}
                    \{ 2 \lambda_j(\SSigma_k^{1/2}\SSigma_{k,0}^{-1}\SSigma_k^{1/2}) -1\} \big]^{-1/2} \label{exp pokazatel}\\
   &\times& \nonumber \exp
                \{ b(k)^{\T} \SSigma_{k,0}^{-1/2}
                        \big[ 2 \SSigma_{k,0}^{-1/2} \SSigma_k^{1/2}
                                \big(2 \SSigma_k^{1/2} \SSigma_{k,0}^{-1}  \SSigma_k^{1/2}
                                -I_{pk}\big)^{-1} \SSigma_k^{1/2} \SSigma_{k,0}^{-1/2}  - I_{pk}\big]
                \SSigma_{k,0}^{-1/2} b(k)
                \} .
\end{eqnarray}
To estimate the obtained expression in terms of the level of noise misspecification~\( \delta \)
 notice that the condition \eqref{multi cond sigma} implies
\begin{equation*}
    \left( \frac{1}{1+\delta} \right)^{pk}
        \le
            \frac{\det \SSigma_k}{\det \SSigma_{k,0}}
        \le
    \left( \frac{1}{1-\delta} \right)^{pk},
\end{equation*}
\begin{equation*}
     \left(\frac{1-\delta}{1+\delta} \right)^{\frac{pk}{2}}
        \le
            \big[ \prod_{j=1}^{pk}
                    \{ 2 \lambda_j(\SSigma_k^{1/2}\SSigma_{k,0}^{-1}\SSigma_k^{1/2}) -1\} \big]^{-1/2}
        \le
     \left( \frac{1+\delta}{1-\delta} \right)^{\frac{pk}{2}}.
\end{equation*}
\begin{equation*}
    \frac{1-\delta}{1+\delta} I_{pk}
        \preceq
            \left( 2 \SSigma_k^{1/2}\SSigma_{k,0}^{-1}\SSigma_k^{1/2} - I_{pk}  \right)^{-1}
        \preceq
    \frac{1+\delta}{1-\delta} I_{pk}.
\end{equation*}
Therefore the quantity in the exponent in \eqref{exp pokazatel} is bounded by:
\begin{eqnarray*}
 &  &  \;\;   \left( 2 \frac{1-\delta}{(1+\delta)^2 } -1 \right)  b(k)^{\T} \SSigma_{k,0}^{-1} b(k)\\
 && \le    b(k)^{\T} \SSigma_{k,0}^{-1/2}
                        \big[ 2 \SSigma_{k,0}^{-1/2} \SSigma_k^{1/2}
                                \big(2 \SSigma_k^{1/2} \SSigma_{k,0}^{-1}  \SSigma_k^{1/2}
                                -I_{pk}\big)^{-1} \SSigma_k^{1/2} \SSigma_{k,0}^{-1/2}  - I_{pk} \big]
                \SSigma_{k,0}^{-1/2} b(k) \\
 & & \le   \left( 2 \frac{1+\delta}{(1-\delta)^2} -1 \right) b(k)^{\T} \SSigma_{k,0}^{-1} b(k).
\end{eqnarray*}
Moreover,
\begin{eqnarray*}
&& \frac{\Delta(k)}{1+\delta}
    =     \frac{1}{1+\delta} b(k)^{\T} \SSigma_k^{-1} b(k) \\
&& \le b(k)^{\T} \SSigma_{k,0}^{-1} b(k)\\
&&  \le \frac{1}{1-\delta} b(k)^{\T} \SSigma_k^{-1} b(k)
    = \frac{\Delta(k)}{1-\delta}.
\end{eqnarray*}
Finally,
\begin{eqnarray}\label{bound for Zk}
  && 
        \left( \frac{1- \delta}{(1+\delta)^3} \right)^{\frac{pk}{2} }
        \exp\left\{ \left(  \frac{2(1-\delta)}{(1+\delta)^2} -1 \right) \frac{\Delta(k)}{1+\delta}\right\}\nn
   &&\le  \EE_{\tta, \Sigma } [Z_k^2]
        \le 
            \left( \frac{1+ \delta}{(1-\delta)^3} \right)^{\frac{pk}{2}}
            \exp\left\{ \left( \frac{2(1+\delta)}{(1-\delta)^2} -1 \right) \frac{\Delta(k)}{1-\delta}\right\}.
\end{eqnarray}

\par In the \emph{case of homogeneous errors} the  expression for \( \log Z_k \) reads as
\begin{eqnarray*}
  \log Z_k  &=&  p k \log \big(\frac{\sigma}{\sigma_0} \big)
                    + \frac{1}{2} \big(\frac{1}{ \sigma^2}-\frac{1}{ \sigma^2_0}\big)
                                            \| \VV_k^{-1/2} (y - \vec \TTa_k) \|^2 \\
            &+& \frac{1}{ \sigma^2_0} b(k)^{\T} \VV_k^{-1} (y-\vec \TTa_k)
                    - \frac{1}{2 \sigma^2_0} b(k)^{\T} \VV_k^{-1} b(k),
\end{eqnarray*}
implying
\begin{equation*}
    \EE_{\tta, \sigma} [Z_k^2] = \left( \frac{\sigma^2}{\sigma_0^2} \right)^{pk}
        \left( \frac{\sigma_0^2}{2 \sigma^2 -\sigma_0^2}  \right)^{\frac{pk}{2}}
        \exp \left\{  \frac{b(k)^{\T} \VV_k^{-1} b(k) }{ 2 \sigma^2 -\sigma_0^2 } \right\}.
\end{equation*}
By Assumption \( (A3)  \)
\begin{eqnarray}\label{bound for Zk homog}
&&
                \left( \frac{1-\delta}{(1+\delta)^3}  \right)^{\frac{pk}{2}}
                \exp \left\{ \frac{\Delta_1(k)}{\sigma^2 (1+\delta)}   \right\}\nn
 \le \EE_{\tta, \sigma} [Z_k^2]
    &\le&  
                \left( \frac{1+\delta}{(1-\delta)^3}  \right)^{\frac{pk}{2}}
                \exp \left\{ \frac{\Delta_1(k)}{\sigma^2 (1-\delta)}   \right\} ,
\end{eqnarray}
where \( p \) is the dimension of the parameter set and \( k \) is the degree of the localization.
\end{proof}
\subsection{Bounds for the bias and variance} \label{proof bias-var}

 Before proceeding with the proof we need to show that the weights \( W^*_{l,\,i}(x) \) defined by \eqref{polynomial weights_multiscaled}
preserve the reproducing polynomials property:
\begin{lemma}\label{Henderson th multiscaled}
Let \( x \in \RR\) be such that Assumptions \( (A1)-(A2) \) hold. Then the weights defined by \eqref{polynomial weights_multiscaled} satisfy
\begin{eqnarray}\label{normalization of polynomial weights multiscaled}
   & & \sum_{i=1}^{n} W^*_{l,\,i}(x) = 1,  \\ \nonumber
   & & \sum_{i=1}^{n} (\Xi - x)^m W^*_{l,\,i}(x) = 0 \; , \;\; m = 1, \ldots, p-1.
\end{eqnarray}
for all \( l=1, \ldots , K \) and any design points \( \{ X_1, \ldots, X_n \} \).
\end{lemma}
\begin{proof}
The assertion can be easily obtained similarly to the proof of Proposition 1.12 from \cite{Tsybakov}.
\end{proof}

\begin{proof} \emph{of Lemma} \ref{nonuniform upper bounds for bias and variance}.
\par
By Lemma \ref{Henderson th multiscaled} and the Taylor theorem with \( \tau_i \) such that the points \( \tau_i \Xi \) are between \( \Xi \) and \( x \), and utilizing Assumption \( (A10) \) we have with \( b_{l,f}(x) =\bb e_1^{\T} \bbpf{l}(x) -  f(x) \):
\begin{eqnarray*}
   |b_{l,f}(x)|   &\le & \frac{1}{(p-1)!}
                        \sum_{i=1}^{n} | f^{(p-1)}(\tau_i \Xi) -f^{(p-1)}(x)|                         |\Xi -x|^{p-1} |W^*_{l,\,i}(x)|\\
   &\le & \frac{L}{(p-1)!}
                        \sum_{i=1}^{n}  |\tau_i \Xi -x|^{\beta -(p-1)} |\Xi -x|^{p-1}|W^*_{l,\,i}(x)|\\
            &\le & \frac{Lh_l^{\beta}}{(p-1)!} \sum_{i=1}^{n}  |W^*_{l,\,i}(x)|.
\end{eqnarray*}
Under the assumptions of the theorem the sum of the polynomial weights can be bounded as follows:
\begin{eqnarray*}
  \sum_{i=1}^{n}  |W^*_{l,\,i}(x)|
        &\le & w_{max}
            \sum_{i=1}^{n} \s^{-2}_i \| \B{l}^{-1} \Psii \| \\
        &\le &
            \kappa(\Sigma) \frac{w_{max}}{\lambda_0 nh_l}  \sum_{i=1}^n \| \Psii \| \,
                        \ind\{ \Xi \in [x-h_l,x+h_l] \}\\
       &\le &
            \kappa(\Sigma)
                \frac{w_{max} \sqrt{e}}{\lambda_0 } a_0 \max\{ 2, \frac{1}{nh_l} \}\\
        &\le &
            \kappa(\Sigma)\frac{2 a_0 w_{max} \sqrt{e} }{\lambda_0 },
\end{eqnarray*}
and the first assertion is justified in view of:
\begin{equation}
    \bar{b}_{k,f}(x) \eqdef \sup_{1\le l \le k } |b_{l,f}(x)|
                            \le \kappa(\Sigma)\frac{2 a_0 w_{max} \sqrt{e} a_0}{\lambda_0 } \frac{Lh_k^{\beta}}{(p-1)!}.
\end{equation}
To bound the variance, just notice that by \eqref{bound for Bk via the smallest eigenvalue in Rd} for any \( \gamma \in \RRp \)
\begin{equation*}
    \gamma^{\T} \B{k}^{-1} \gamma \le \frac{\s^2_{max} }{nh_k \Lambda_0} \| \gamma \|^2.
\end{equation*}
Then under Condition \( (A3)\) by \eqref{Vk bound} for the variance term we have:
\begin{eqnarray*}
  \s^2_k(x) &=& \bb e_1^{\T} \Var\mmle{k} \, \bb e_1 \\
  &\le & (1+\delta) \bb e_1^{\T} \B{k}^{-1} \bb e_1 \\
  &\le & (1+\delta) \frac{\s^2_{max} }{nh_k \Lambda_0}.
\end{eqnarray*}

\end{proof}


\begin{thebibliography}{10}

{

\bibitem{Akaike}
Akaike, H. (1973).
Information theory and an extension of the maximum
likelihood principle. \textit{ Second International Symposium on Information
Theory (Tsahkadsor, 1971), Akad\'emiai Kiad\'o, Budapest} 267--281.


\bibitem{Arlot}
Arlot, S. (2009). Model selection by resampling penalization. \emph{Electron. J. Stat.} \textbf{3} 557--624.

\bibitem{ArlotMassart}
Arlot, S. and Massart, P. (2009). Data-driven calibration of penalties for least squares regression.
\emph{J. Mach. Learn. Res.} \textbf{10}(Feb) 245--279.

\bibitem{BaraudGiraudHuet}
Baraud, Y., Giraud, C. and Huet, S. (2009). Gaussian model selection with an unknown variance.
\textit{ Ann. Statist.} \textbf{37:2} 630--672.

\bibitem{BirgeMassart}
Birg\'e, L. and Massart, P. (2001). Gaussian model selection.
\emph{Journal of the European Mathematical Society} \textbf{3:3} 203--268.

\bibitem{Brillinger}
Brillinger, D. R. (1977). Discussion of Stone (1977). \textit{ Ann. Statist.} \textbf{5:4} 622--623.

\bibitem{Brua}
Brua, J.-Y. (2009). Asymptotic efficient estimators for non-parametric
heteroscedastic model. \emph{Statistical Methodology} \textbf{6:1} 47--60.

\bibitem{Cleveland79}
Cleveland, W. S. (1979). Robust locally weighted regression and smoothing scatterplots. \emph{J. Amer. Statist. Assoc.} \textbf{74:368} 829--836.

\bibitem{DalalyanSalmon}
Dalalyan A. S. and Salmon J. (2011). Sharp Oracle Inequalities for Aggregation of Affine Estimators. Preprint arXiv:1104.3969v2.

\bibitem{Donoho and Johnstone}
Donoho, D. L. and Johnstone, I. M. (1994).
Ideal spatial adaptation by wavelet shrinkage. \textit{Biometrica} \textbf{81} 425--455

\bibitem{Fan and Gijbels book}
Fan, J. and Gijbels, I. (1996). \textit{Local Polynomial Modelling and Its Applications.}
Monographs on Statistics and Applied Probability, 66. Chapman and Hall,
London.

\bibitem{FanZhangZhang}
Fan, J., Zhang, C. and Zhang, J. (2001).
Generalized likelihood ratio statistics and Wilks phenomenon. \textit{ Ann. Statist.}  \textbf{29:1} 153--193.

\bibitem{GaltchoukPergam2009}
Galtchouk, L. and Pergamenshchikov, S. (2009).
Adaptive asymptotically efficient estimation in heteroscedastic nonparametric regression
\emph{Journal of the Korean Statistical Society} \textbf{38:4} 305--322.


\bibitem{GaltchoukPergam2010a}
Galtchouk, L. and Pergamenshchikov, S. (2010).  Adaptive asymptotically efficient estimation
in heteroscedastic nonparametric regression. arXiv:1002.1537v1

\bibitem{GaltchoukPergam2010b}
Galtchouk, L. and Pergamenshchikov, S. (2010).  Sharp non-asymptotic oracle inequalities for
nonparametric heteroscedastic regression models. arXiv:1002.1538v1



\bibitem{Goldenshluger and Nemirovski}
Goldenshluger, A. and Nemirovski, A. (1994). On spatial adaptive estimation of nonparametric regression.
\textit{Research report, Technion-Israel Inst. Technology, Haifa, Israel}.

\bibitem{Golubev_Spokoiny}
Golubev, Y. and Spokoiny, V. (2009). Exponential bounds for minimum contrast estimators. \emph{Electron. J. Stat.} \textbf{3} 712--746.

\bibitem{EfroimovichPinsker96}
 Efroimovich, S. and Pinsker, M.(1996).
 Sharp-optimal and adaptive estimation for heteroscedastic nonparametric regression.
 \emph{Statistica Sinica} \textbf{6} 925--942.

\bibitem{Efroimovich2007}
Efroimovich, S. (2007).
Sequential design and estimation in heteroscedastic nonparametric regression. \emph{Sequential Analysis} \textbf{26} 3--25.

\bibitem{Katk1979}
Katkovnik, V. Ja. (1979).
Linear and nonlinear methods of nonparametric regression analysis. (Russian)  \textit{Soviet Automat. Control} \textbf{5} 35--46, 93.

\bibitem{Katk1983}
Katkovnik, V. Ja. (1983).
Convergence of linear and nonlinear nonparametric estimates of ``kernel'' type.  \textit{Automat. Remote Control}  \textbf{44:4} 495--506;  translated from  \textit{Avtomat. i Telemekh.}  1983 \textbf{4} 108--120 (Russian).


\bibitem{Katk1985}
Katkovnik, V. Ja. (1985).
\textit{Nonparametric Identification and Data Smoothing: Local Approximation Approach}. Nauka, Moscow (Russian).

\bibitem{KatkEA2006}
 Katkovnik, V., Egiazarian, K. and Astola, J. (2006).
 \textit{Local Approximation Techniques in Signal and Image Processing}. Bellingham, WA: SPIE Press.

\bibitem{KatkSpok}
Katkovnik, V. and Spokoiny, V. (2008). Spatially adaptive estimation via fitted local likelihood techniques. \textit{IEEE Trans. Signal Process.}, \textbf{56}, No.3, 873--886.

\bibitem{KLP2001}
 Kerkyacharian, G., Lepski, O. and Picard, D. (2001) Nonlinear estimation in anisotropic multi-index denoising.
 \emph{Probab. Theory Related Fields} \textbf{121:2} 137--170.

 \bibitem{KLP2007}
 Kerkyacharian, G., Lepski, O. and Picard, D. (2007) Nonlinear estimation in anisotropic multiindex denoising.
 Sparse case. \emph{Teor. Veroyatn. Primen.} \textbf{52:1} 150--171; translation in
 \emph{Theory Probab. Appl.} (2008) \textbf{52:1} 58--77.

\bibitem{KullbackLeibler}
Kullback, S. and Leibler, R. A. (1951).
On information and sufficiency.
\textit{Ann. Math. Statistics} \textbf{22} 79--86.


 \bibitem{Lep1990}
 Lepskii, O. V. (1990).
 A problem of adaptive estimation in Gaussian white noise.  (Russian) \textit{Teor. Veroyatnost. i Primenen.}  \textbf{35:3} 459--470;  translation in \textit{Theory Probab. Appl.} \textbf{35:3} 454--466.


 \bibitem{Lep1992}
 Lepskii, O. V. (1992).
 Asymptotic minimax adaptive estimation. II. Schemes without optimal adaptation. Adaptive estimates. (Russian) \textit{Teor. Veroyatnost. i Primenen.} \textbf{37:3} 468--481; translation in
\textit{Theory Probab. Appl.} \textbf{37:3} 433--448.

\bibitem{LepMamSpok97}
Lepski, O. V., Mammen, E. and Spokoiny, V.G. (1997).
Optimal spatial adaptation to inhomogeneous smoothness: an approach based on kernel estimates with variable bandwidth selectors. \textit{Ann. Stat.} \textbf{25:3} 929--947.

\bibitem{LepSpok97}
Lepski, O. V. and Spokoiny, V.G. (1997).
Optimal pointwise adaptive methods in nonparametric estimation. \textit{Ann. Stat.} \textbf{25:6} 2512--2546.

\bibitem{Loader}
Loader, C. (1999). \textit{Local Regression and Likelihood.} Statistics and Computing. Springer-Verlag, New York.

\bibitem{Massart}
Massart, P. (2003). \emph{Concentration Inequalities and Model Selection} (2007). Ecole d'\'et\'e de Probabilit\'es de Saint-Flour . Lecture Notes in Mathematics 1896, Springer Berlin/Heidelberg.

\bibitem{RuppertWand}
Ruppert, D. and Wand, M. P. (1994). Multivariate locally weighted least squares regression. \textit{Ann. Stat.} \textbf{22:}3 1346--1370.

\bibitem{Saumard}
Saumard, A. (2010). Optimal upper and lower bounds for the true and empirical excess risks in heteroscedastic least-squares regression, Preprint hal-00512304, v1.

\bibitem{Spokoiny variance}
 Spokoiny, V. (2002).
 Variance estimation for high-dimensional regression models. \textit{J. Multivariate Anal.}
 \textbf{82} 111--133.

\bibitem{SV}
 Spokoiny, V. and Vial, C. (2009). Parameter tuning in pointwise adaptation using a propagation approach. \textit{Ann. Statist.} \textbf{37:5B} 2783--2807.

\bibitem{Stone77}
Stone, C. J. (1977). Consistent nonparametric regression. With discussion and a reply by the author. \emph{Ann. Statist.} \textbf{5:4} 595--645.

\bibitem{Stone80}
Stone, C. J. (1980). Optimal rates of convergence for nonparametric estimators. \emph{Ann. Statist.} \textbf{8:6} 1348--1360.

\bibitem{TibshiraniHastie}
Tibshirani, R. and Hastie, T. (1987). Local likelihood estimation. \textit{J. Amer. Statist. Assoc.} \textbf{82:398} 559--567.

\bibitem{Tsybakov86}
Tsybakov, A. B. (1986). Robust reconstruction of functions by a local approximation method. (Russian) \emph{Problemy Peredachi Informatsii} \textbf{22:2} 69--84 (\emph{Problems of Information Transmission}, 1986 \textbf{22:2} 133--146 ).

\bibitem{Tsybakov}
Tsybakov, A. B. (2009). \textit{Introduction to Nonparametric Estimation}.
Springer Series in Statistics. Springer-Verlag, New York.

\bibitem{White}
White, H. (1982).
Maximum likelihood estimation of misspecified models. \textit{Econometrica}  \textbf{50:1} 1--25.
}
\end{thebibliography}
\end{document}